\documentclass[a4paper]{article}
\usepackage{academicPaper}

\newcommand{\tor}{\text{\textbf{-tor}}}
\DeclareMathOperator{\tra}{Tra}
\DeclareMathOperator{\hol}{Hol}

\DeclareMathOperator{\Hom}{Hom}
\DeclareMathOperator{\Trans}{\mathbf{Trans}}
\DeclareMathOperator{\Bun}{\mathbf{Bun}}
\renewcommand{\hat}{\widehat}
\newcommand{\thin}{\text{thin}}

\begin{document}

\title{\Large\uppercase{\textbf{Parallel 2-transport and 2-group torsors}}}
\author{\textsc{W.H. Voorhaar}\\ Université de Genève, \href{mailto:wh.voorhaar@gmail.com}{wh.voorhaar@gmail.com}}
\date{}

\maketitle

\begin{abstract}
	We provide a new perspective on parallel 2-transport and principal 2-group bundles with 2-connection. We define parallel 2-transport as a 2-functor from the thin fundamental 2-groupoid to the 2-category of 2-group torsors. The definition of the 2-category of 2-group torsors is new, and we develop the tools necessary for computations in this 2-category. We prove a version of the non-Abelian Stokes Theorem and the Ambrose-Singer Theorem for 2-transport. This definition  motivated by the fact that principal $G$-bundles with connection are equivalent to functors from the thin fundamental groupoid to the category of $G$-torsors. In the same lines we deduce a notion of principal 2-bundle with 2-connection, and show it is equivalent to our notion 2-transport functors. This gives a stricter notion than appears in the literature, which is more concrete. It allows for computations of 2-holonomy which will be exploited in a companion paper to define Wilson surface observables. Furthermore this notion can be generalized to a concrete but strict notion of $n$-transport for arbitrary $n$. 
\end{abstract}

\tableofcontents

\newpage

\section*{Acknowledgments}
This research was supported by the National Center for Competence in Research SwissMAP of the Swiss National Science Foundation. I would like to express my gratitude towards Eugene Lerman for his invaluable help and guidance, particularly at the early stages of this project. A significant amount of content of this paper is derived from my master thesis. I would like to thank Anton Alekseev for his role as adviser of my master thesis, and for his guidance both during and after the writing of my master thesis. Finally I would like to thank NCCR SwissMAP for hosting me in Geneva during time of writing my master thesis.

\section{Introduction}

There is a great variety of notions of parallel 2-transport for principal 2-bundles with 2-connection, consider for instance \cite{Wal17,SchWal13,CattaneoCottaRamusinoRinaldi,MartinsPicken1,MartinsPicken2,SW11,SW13,Wang}. The majority of these papers define 2-transport locally and then use a gluing construction to extend to the global case. Instead we will develop a notion that is inherently global, based on considering 2-functors $\Pi^\thin_2(M)\to \mc G\tor$. This is inspired by \cite{CollierLermanWolbert} which proves an equivalence between functors $\Pi^\thin_1(M)\to G\tor$ and principal $G$-bundles with connection, and \cite{SW11} which develops theory for trivial principal $\mc G$-2 bundles with 2-connection by considering functors $\Pi^\thin_2(M)\to B\mc G$.

This definition arose from a need to have a notion of 2-transport which can deal with non-trivial 2-bundles, but which is also more concrete than appears elsewhere in the literature. This in turn is based on the observation that in practice all the principal 2-bundles we encountered were based on constructions out of ordinary principal bundles. Such principal 2-bundles tend to have the same local data as ordinary principal bundles, and are stricter than most notions of 2-bundles appearing in the literature. Hence we develop a notion of 2-transport that is more suited to the structures we encounter in practice. 

The main application of this theory will be to develop a mathematical notion of Wilson surface observables, to appear in a separate paper. That is, a (higher) gauge invariant quantity associated to a map $S\to M$ from a surface to a manifold with a principal (2-)bundle. The idea there is to interpret the results in \cite{AlekseevChekeresMnev} in the language of higher gauge theory, and to extend the notion obtained in this way to different 2-groups. 

The theory of 2-gauge theory can also be seen as a bridge towards understanding $n$-gauge theory. Indeed, armed with a higher version of the non-Abelian Stokes theorem, it becomes fairly straightforward to define principal $n$-bundles with $n$-connections and their associated parallel $n$-transport, at least for strict $n$-groups. This will be explained in a separate paper. 

\section{General theory}

\subsection{Lie 2-groups}
A Lie 2-group is a 2-categorical generalization of a Lie group. There are several equivalent ways of defining 2-groups. We will be interested in \textit{strict} Lie 2-groups only.

\begin{definition}
 Abstractly, a Lie 2-group is a category internal to \textbf{LieGrp}, the category of Lie groups. Or a (strict, smooth) 2-category with one object and all morphisms invertible. More concretely it is a Lie groupoid $\mc G_1\rightrightarrows \mc G_0$ with a group multiplication $\mc M\colon\mc G_i\times \mc G_i\to \mc G_i$, which is compatible with the groupoid structure (i.e. group multiplication and inversion are functors).
\end{definition}

\noindent Yet even more concretely, a Lie 2-group is equivalent to a smooth crossed module.

\begin{definition}
	A smooth crossed module $(G,H,t,\alpha)$ consists of two Lie groups $G$, $H$, together with a Lie group homomorphism $t\colon H\to G$ (the `target map') and an action $\alpha\colon G\to \Aut(H)$ satisfying the following two axioms:
	\begin{itemize}
		\item $t (\alpha_gh) = gt(h)g^{-1}$ for all $g\in G$, $h\in H$
		\item $\alpha_{t(h)}h' = hh'h^{-1}$ for all $h,h'\in H$.
	\end{itemize}
\end{definition}

	\begin{proposition}
		There is a 1:1 equivalence between smooth crossed modules and Lie 2-groups.
	\end{proposition}
	\proof[Proof (sketch):]
	A smooth crossed module gives a Lie 2-group by setting $\mc G_1 = G\ltimes H$ and $\mc G_0 = G$ with source $(g,h)\mapsto g$, target $(g,h)\mapsto t(h)g$, multiplication $(g,h)(g',h') = (gg',h\alpha_g(h'))$ and composition $(t(h)g,h')\circ (g,h) = (g,h'h)$. Conversely given a Lie 2-group $s,t\colon \mc G_1\rightrightarrows \mc G_0$ we get a crossed module by setting $G=\mc G_0$, $H = \ker s$, $t=t|_{\ker s}$, and $\alpha_g(h) = \id_g\cdot h\cdot \id_g^{-1}$. Consult \cite{BrownSpencer} for more details.\qed
	
	From now one we will treat smooth crossed modules and Lie 2-groups without distinction, and the words 'smooth' and 'Lie' are sometimes omitted. Moreover we fix a smooth crossed module $(G,H,t,\alpha)$ which is equivalent to the 2-group $G\ltimes H\rightrightarrows G$ for the rest of this paper.

	\begin{proposition}
		For any crossed module $\ker t\subset H$ is a central subgroup.
	\end{proposition}
	\proof Let $h\in \ker t$ then $hh'h^{-1} = \alpha_{t(h)}h'=h'$ for any $h'\in H$.\qed

	\noindent Conversely any central extension of Lie groups
	\[
		1\to Z\to H\to G\to 1
	\]
	defines a crossed module. Note that this is equivalent to $t$ being surjective. This gives a very large class of examples of 2-groups, and many 2-groups that arise in practice are of this form.
	
	The infinitesimal counterpart of a Lie 2-group is a Lie 2-algebra. We will describe Lie 2-algebras from the point of view of crossed modules.
	\begin{definition}
		The \textit{Lie 2-algebra} or \textit{infinitesimal crossed module} associated to a crossed module $(G,H,t,\alpha)$ is the tuple $(\mf g,\mf h,t_*,\alpha_*)$ where $\mf g$, $\mf h$ are the Lie algebras of $G$, $H$ respectively, $t_*:\mf h\to \mf g$ is the derivative of $t:H\to G$, and $\alpha_*:\mf g\times \mf h\to \mf h$ is the bilinear map obtained by differentiating $\alpha:G\times H\to H$ in both arguments.
	\end{definition}
	
	To an infinitesimal crossed module we also associate the Lie algebra $\mf g\ltimes \mf h$, which as a vector space is $\mf g\times \mf h$ but with Lie bracket:
	\begin{align}\label{eq:infcrossedmodulebracket}
		&[X,Y]_{\mf g\ltimes \mf h} = [X,Y]_{\mf g}
		&& [X,\xi]_{\mf g\ltimes \mf h} = \alpha_*(X,\xi)
		&&[\xi,\eta]_{\mf g\ltimes \mf h} = [\xi,\eta]_{\mf h} &X,Y\in \mf g,\,\,\xi,\eta\in\mf h.
	\end{align} 
	We note that $s_*(X+\xi)=X$, $t_*(X+\xi)=X+t_*\xi$ are source and target maps for a groupoid $\mf g\ltimes \mf h\rightrightarrows \mf g$.

	\subsection{$\mc G$-2-torsors}
	
	A $G$-torsor, for $G$ a Lie group is a manifold with a free and transitive $G$-action, or equivalently, a principal $G$-bundle over a point. Similarly, $\mc G$-2-torsors are Lie groupoids with a free and transitive $\mc G$-action. These $\mc G$-2-torsors will form the fibers of principal 2-bundles, as well as the target 2-category of parallel 2-transport. First let us recall the definition of a 2-group action on a Lie groupoid.
	
	\begin{definition}
		A (right) Lie 2-group action of $\mc G=(G\ltimes H\rightrightarrows G)$ on a Lie groupoid $\mc X_1\rightrightarrows \mc X_0$ is a functor $R\colon  \mc X\times \mc G\to \mc X$ such that the following diagram commutes (on the nose):
	\[
		\begin{tikzcd}[column sep=4em, row sep =4em]
		 \mc X\times\mc G\times \mc G\ar["  \id\times M"]{r} \ar["R \times \id"']{d} & \mc X\times \mc G\ar["R"]{d}\\
			\mc X\times \mc G \ar["R"']{r} & \mc X	
		\end{tikzcd}
		\]
	where $M\colon \mc G\times \mc G \to \mc G$ is the 2-group multiplication. That is, there is a $G\ltimes H$-action on $\mc X_1$ and a $G$ action on $\mc X_0$ which are compatible with the groupoid structure.
	\end{definition}

	\begin{definition}
		Let $\mc G$ be a 2-group, then a \textit{$\mc G$-2-torsor} is a groupoid $\mc X_1\rightrightarrows\mc X_0$ together with a free and transitive $\mc G$ action $R\colon \mc X\times \mc G\to \mc X$. That is, the $G\ltimes H$ and $G$ actions are both free and transitive, or in other words $\mc X_1$ is a $G\ltimes H$-torsor and $\mc X_0$ is a $G$-torsor. Equivalently the functor
		\[
			(\pi_1\times R)\colon \mc X\times G\to \mc X\times \mc X
		\]
		is invertible (on the nose), where $\pi_1$ is the projection onto the first factor. That is, for $X,Y\in \mc X_i$ there is a unique $Y:X\in \mc G_i$ such that $X\cdot(Y:X)=Y$.
	\end{definition}
	
	This definition has a number of useful consequences. First we derive two useful formulas for division and composition. Note that the map 
	\begin{equation}
	\mc X\times \mc X\to \mc G,\qquad (X,Y)\mapsto X:Y
	\end{equation}
	is functorial. In particular this means that for any $X,X',Y,Y'\in \mc X$ such that $Y\circ X$ and $Y'\circ X'$ are defined we have
	\begin{equation}\label{eq:DivisionFunctorial}
		(Y:Y')\circ (X:X') = (Y\circ X):(Y'\circ X').
	\end{equation}
	Functoriality of the $\mc G$ action implies that for any $X,Y\in \mc X$ and $g,h\in \mc G$ such that $X\circ Y$ and $g\circ h$ are defined we have
	\begin{equation}\label{eq:CompositionEquivariant}
		(X\cdot g)\circ (Y\cdot h) = (X\circ Y)\cdot (g\circ h).
	\end{equation}
	
	Next we consider the morphisms of $\mc G$-2-torsors, and elaborate on the properties of these morphisms.

	\begin{definition}
		$\mc G$-2-torsors form a (strict) 2-category $\mc G\tor$. The 1-morphisms are given by smooth equivariant functors, and the 2-morphisms by smooth natural transformations. 
	\end{definition}
	
		\begin{lemma}\label{lem:equivfunctor}
			A 1-morphism of $(G\ltimes H\rightrightarrows G)$-2-torsors $F\colon \mc X\to \mc Y$ is completely determined by the induced morphism $F_0\colon \mc X_0\to \mc Y_0$ of $G$-torsors. Furthermore any morphism $F_0\colon \mc X_0\to \mc Y_0$ extends uniquely to a 1-morphism  $F\colon \mc X\to \mc Y$.
		\end{lemma}
	
		\proof Let $F\colon \mc X\to \mc Y$ be an equivariant functor and let $X\colon p\to q$. Then using functoriality and equivariance we have
		\begin{equation}
			F(X) = F\left(\id_p\cdot\, (X:\id_p)\right) = F(\id_p)\cdot (X:\id_p)=\id_{F_0(p)}\cdot (X:\id_p),
		\end{equation}
		thus $F$ is completely determined by $F_0\colon \mc X_0\to \mc Y_0$. On the other hand suppose $F_0\colon \mc X_0\to \mc Y_0$ is an equivariant map, then we claim
		\begin{equation}
			F(X) = \id_{F_0(p)}\cdot (X:\id_p)
		\end{equation}
		defines an equivariant functor $F\colon \mc X\to \mc Y$. Let $X\colon p\to q$ and $Y\colon q\to r$, then we need to show that
		\begin{equation}\label{eq:FFunctorial}
			F(Y\circ X) = F(Y)\circ F(X).
		\end{equation}
		We now expand the left hand side of \eqref{eq:FFunctorial} and divide by $\id_{F(p)}$:
		\begin{align}
				F(Y\circ X):\id_{F(p)} = (Y\circ X):\id_p.
		\end{align}
		On the other hand, using \eqref{eq:DivisionFunctorial} twice we compute,
		\begin{align}
			(F(Y)\circ F(X)) : \id_{F(p)} &= (\id_{F(q)}\cdot Y:\id_q)\circ (\id_{F(p)}\cdot (X:\id_p)):\id_{F(p)}\\
			&=(\id_F(q)\cdot (Y:\id_q):\id_{F(p)})\circ (X:\id_p)\\
			&=((\id_{F(q)}:\id_{F(p)})\cdot (Y:\id_q))\circ (X:\id_p)\\
			&=((\id_{q}:\id_{p})\cdot (Y:\id_q))\circ (X:\id_p)\\
			&=(Y:\id_p)\circ(X:\id_p) = (Y\circ X):\id_p,
		\end{align}
		which shows that \eqref{eq:FFunctorial} holds, and hence that $F$ is a functor. Equivariance is per definition.\qed
		
			\begin{lemma}\label{lem:g2tor2morph}
				One can equivalently define a 2-morphism $F\Rightarrow F'$ for $F,F'\colon \mc X\to \mc Y$ as a map $\eta\colon \mc X_0\to \mc Y_1$ such that:
				\begin{enumerate}
					\item $\eta(p):F(p)\to F'(p)$
					\item $\eta(p\cdot g) = \eta(p)\cdot \id_g$
				\end{enumerate}
			\end{lemma}
			\proof The first identity is per definition of a natural transformation. We will show the second identity is equivalent to the property that for any $X\colon p\to q$ in $\mc X$ the following diagram commutes:
			\begin{equation}
				\begin{tikzcd}[column sep=4em, row sep =4em]
					F(p) \ar["F(X)"]{r} \ar["\eta(p)"']{d} & F(q) \ar["\eta(q)"]{d}\\
					F'(p) \ar["F'(X)"']{r} & F'(q)
				\end{tikzcd}
			\end{equation}
			Writing this out as a formula we obtain
			\begin{equation}\label{eq:NaturalTransformCommute}
				\eta(q)\circ (\id_{F(p)}\cdot (X:\id_p)) = (\id_{F'(p)}\cdot (X:\id_p))\circ \eta(p).
			\end{equation}
			We apply equivariance of composition \eqref{eq:CompositionEquivariant} to the left hand side of \eqref{eq:NaturalTransformCommute}:
			\begin{align}
					\eta(q)\circ (\id_{F(p)}\cdot (X:\id_p)) &= \left(\left[\eta(q)\cdot \id_{F(p):F(q)}\right]\cdot \id_{F(q):F(p)}\right)\circ (\id_{F(p)}\cdot (X:\id_p))\\
					&=\left(\eta(q)\cdot \id_{F(p):F(q)})\right)\cdot  X:\id_p.
			\end{align}
			We can apply the same trick to the right hand side of \eqref{eq:NaturalTransformCommute}:
			\begin{align}
				 (\id_{F'(p)}\cdot (X:\id_p))\circ \eta(p) =  (\id_{F'(p)}\cdot (X:\id_p))\circ \eta(p)\cdot \id_1=\eta(p)\cdot (X:\id_p).
			\end{align}
			Thus we obtain that \eqref{eq:NaturalTransformCommute} is equivalent to
			\begin{equation}
				\eta(p) = \eta(q)\cdot \id_{F(p):F(q)}=\eta(q)\cdot \id_{p:q},
			\end{equation}
			which is precisely the identity we wished to prove. \qed
			
			\begin{lemma}\label{lem:g2tor2morphH}
				The map $\eta\colon \mc X_0\to \mc Y_1$ of lemma \ref{lem:g2tor2morph} is equivalent to the map $\eta_H\colon \mc X_0\to H$ defined by
				\begin{equation}\label{eq:etaH}
					\eta_H(p) = \eta(p):\id_{F(p)}.
				\end{equation}
				In fact a 2-morphism $F\rightarrow F'$ is equivalent to a map $\eta_H\colon \mc X_0\to H$ such that 
				\begin{enumerate}
					\item $t(\eta_H(p)) = F'(p):F(p)$
					\item $\eta_H(p\cdot g) =  \alpha_{g^{-1}}\eta_H(p)$
				\end{enumerate}
			\end{lemma}
			\proof Equation \eqref{eq:etaH} together with lemma \ref{lem:g2tor2morph} provide the equivalence. The first property in this lemma is immediate, the second is a computation:
			\begin{align*}
				\eta_H(p\cdot g) &= \eta(p\cdot g):\id_{F(p\cdot g)}\\
				&= (\eta(p)\cdot 1_g):(\id_{F(p)}\cdot 1_g) \\
				&= 1_{g^{-1}}\cdot (\eta(p):\id_{F(p)})\cdot 1_g \\
				&= \alpha_{g^{-1}}\eta_H(p).\qedhere
			\end{align*}
			
	\begin{lemma}\label{lem:etaHcomposition}
		Let $\bullet$ denote vertical composition, and $\circ$ horizontal composition of 2-morphisms of $\mc G$-torsors. Then horizontal and vertical composition in terms of $\eta_H$ is given by:
		\begin{align}
			\hspace{-2em}\begin{tikzcd}[column sep=6em, row sep =6em,ampersand replacement=\&]
				\mc X \ar[bend left = 60,"F"{name = A}]{r} \ar["F'"{description, name=B}]{r} \ar[bend right = 60,"F''"'{name=C}]{r} \ar[Rightarrow, from=A, to=B,shorten <=0.1cm, "\eta"{pos=0.6}] \ar[Rightarrow, from=B, to=C,shorten >=0.1cm, "\eta'"{pos=0.3}]\& \mc Y
			\end{tikzcd}\qquad &\hspace{4em}
			(\eta'\bullet\eta)_H(p) = \eta_H(p)\cdot \eta'_H(p)\\		
			\hspace{-2em}\begin{tikzcd}[column sep =4em, row sep =4em,ampersand replacement=\&]
				\mc X \ar[bend left=50,"F_1"{name=A}]{r} \ar[bend right=50,"F_1'"'{name=B}]{r} \ar[Rightarrow, from=A, to=B,"\eta_1",shorten >=1mm,shorten <=1mm] \& 
				\mc Y \ar[bend left=50,"F_2"{name=C}]{r} \ar[bend right=50,"F_2'"'{name=D}]{r} \ar[Rightarrow, from=C, to=D,"\eta_2",shorten >=1mm,shorten <=1mm] \&
				\mc Z
			\end{tikzcd}&\hspace{4em} 
			\begin{matrix}
			(\eta_2\circ \eta_1)_H(p) \hspace{-0.8em}&= \eta_{1,H}(p)\cdot \eta_{2,H}(F_1'(p))\\
						&=\eta_{2,H}(F_1(p))\cdot\eta_{1,H}(p)
			\end{matrix}
		\end{align}
	\end{lemma}
	
	\proof We compute the first identity:
	\begin{align*}
		(\eta'\bullet \eta)_H(p) &= \left(\eta'(p)\circ \eta(p)\right):\id_{F(p)}=\eta'(p):\id_{F(p)}\circ\, \eta_H(p)\\
		 &= \id_{F'(p):F(p)}\eta'_H(p)\circ \eta_H(p) = \alpha_{t(\eta_H(p))}(\eta'_H(p))\circ \eta_H(p)\\
		 & = \eta_H(p)\eta'_H(p).
	\end{align*}
	For the second identity recall first the definition of horizontal composition of natural transformations:
	\begin{equation}
		(\eta_2\circ\eta_1)(p) = \eta_2(F_1'(p))\circ F_2(\eta_1(p)) = F_2'(\eta_1(p))\circ \eta_2(F_1(p))
	\end{equation}
	Then we compute
	\begin{align*}
		(\eta_2\circ\eta_1)_H(p)&=\big[\eta_2(F_1'(p))\circ F_2(\eta_1(p))\big]:\id_{F_2F_1(p)}\\
		& = \big[\eta_2(F_1'(p)):\id_{F_2F_1(p)}\big]\circ \big[F_2(\eta_1(p)):\id_{F_2F_1(p)}\big]\\
		& = \big[\eta_2(F_1(p))\cdot \id_{F_1'(p):F_1(p)}:\id_{F_2F_1(p)}\big]\circ \big[\id_{F_2F_1(p)}\cdot \eta_{1,H}(p):\id_{F_2F_1(p)}\big]\\
		&= \eta_{2,H}(F_1(p))\cdot \id_{t(\eta_{1,H}(p))}\circ \eta_{1,H}(p)\\
		&= \eta_{2,H}(F_1(p))\cdot \eta_{1,H}(p).
	\end{align*}
	We can also transform this identity, to complete the proof:
	\begin{align*}
		\eta_{2,H}(F_1(p))\cdot \eta_{1,H}(p) &= \eta_{2,H}(F'_1(p)\cdot F_1(p):F_1'(p))\cdot \eta_{1,H}(p)
		\\ &= \alpha(t(\eta_{1,H}(p)),\eta_{2,H}(F'_1(p)))\eta_{1,H}(p)
		\\&=\eta_{1,H}(p)\eta_{2,H}(F'_1(p)).\qedhere
	\end{align*}
	
	\subsection{The thin fundamental 2-groupoid}
	Given a manifold $M$, we can consider its fundamental groupoid $\Pi_1(M)$. This is the groupoid with objects points in $M$ and arrows homotopy classes of paths. Instead we will be interested in the \textit{thin} fundamental groupoid $\Pi_1^\thin(M)$, where we take paths up to \textit{thin homotopy} instead. 
	
	\begin{definition}
		Two paths $\gamma_0,\gamma_1$ are \textit{thinly homotopic} if there is a (smooth) homotopy $\Sigma\colon I^2\to M$, $\Sigma\colon \gamma_0\Rightarrow\gamma_1$ such that the rank of the differential $ D\Sigma\colon TI^2\to TM$ is at most 1 at every point.
	\end{definition}
	In other words a thin homotopy $\Sigma$ `doesn't sweep out any area' and changes $\gamma_i$ only in the direction of $\gamma_i$. In this sense it is essentially a notion of reparametrization of paths.	Any path is thinly homotopic to a path which is locally constant near its endpoints (`sitting instants'). Using this we can define concatenation of thin homotopy classes of paths. 
	
	\begin{definition}
		The \textit{thin fundamental groupoid} $\Pi_1^\thin(M)$ is the groupoid with objects points in $M$ and arrows thin homotopy classes of paths. Composition of arrows is concatenation of paths, and inversion is reversal of paths.
	\end{definition}
	
	This can be upgraded to a \textit{thin fundamental 2-groupoid} by adding bigons (i.e. path homotopies) as 2-morphisms. A bigon homotopy $h\colon \Sigma_0\Rrightarrow\Sigma_1$ between bigons $\Sigma_i\colon \gamma_0\Rightarrow\gamma_1$ is a map $h\colon I^3\to M$ such that $h(i,s,t)=\Sigma_i$, and $h(u,s,t)$ is a bigon $\gamma_0\Rightarrow\gamma_1$ for each $u\in I$. A homotopy of bigons $h\colon \Sigma_0\Rrightarrow\Sigma_1$ is \textit{thin} if the rank of the differential $Dh\colon TI^3\to TM$ is at most 2 at any point. Bigons can be horizontally and vertically composed, and this leads to a 2-groupoid. Refer to \cite{SW11} for a more detailed definition.
	
\begin{figure}[!htb]\centering
			\includegraphics[width = 0.7\textwidth]{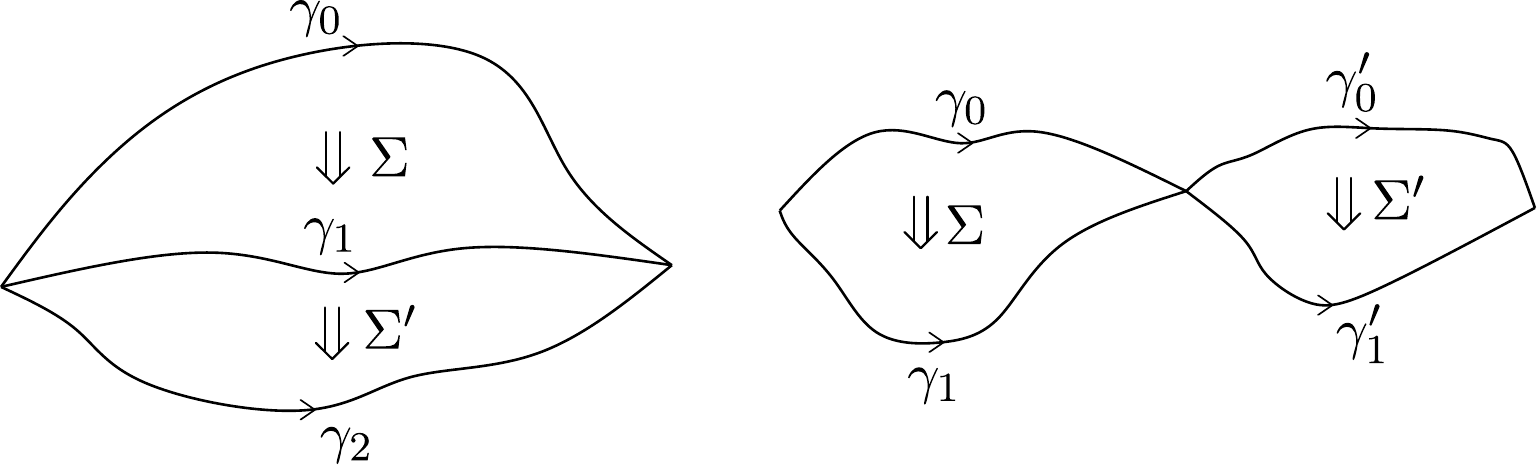}
			\caption{The left shows vertical composition of bigons, and the right side horizontal composition of bigons.\label{fig:bigonComposition}}
		\end{figure}

	\begin{definition}
		The \textit{thin fundamental 2-groupoid} $\Pi_2^\thin(M)$ is the 2-groupoid with objects points in $M$, 1-morphisms thin homotopy classes of paths, and 2-morphisms thin homotopy classes of bigons. Horizontal and vertical composition of bigons is illustrated in figure \ref{fig:bigonComposition}
	\end{definition}
	
	\subsection{Transport 2-functors}
Collier, Lerman and Wolbert \cite{CollierLermanWolbert} showed that smooth functors $\Pi_1^{\thin}(M)\to G\tor$ are the same thing as principal bundles with connection (in fact they established an equivalence of stacks). This motivates us to study smooth 2-functors $\Pi_2^\thin(M)\to \mc G\tor$ as these should give an appropriate notion of principal 2-bundle with 2-connection. So far we have defined the source and target 2-categories of such functors, but it is not immediately obvious what it means for such 2-functors to be smooth. There is a natural diffeological structure on $\Pi_2^\thin(M)$, but for $\mc G\tor$ there is no such structure we can use. Schreiber and Waldorf \cite{SW13} defined smoothness by introducing a notion of locally trivializable 2-functors and smooth descent. We will instead present a generalization to the smoothness notion used in \cite{CollierLermanWolbert}.

	\begin{definition}\label{def:transport2functor}
		A $\mc G$-transport 2-functor is a 2-functor $\ms F\colon \Pi_2^\thin(M)\to \mc G\tor$ satisfying the following smoothness conditions:
	\end{definition} 
		\begin{itemize}
			\item Set $\Pi_1^\thin(M,x)$ to be the group of thin homotopy classes of loops at $x\in M$. Then $\ms F$ induces a homomorphism $\Pi_1^\thin(M,x)\to \Aut(\ms F(x))$. The space on the left is diffeological, and the space on the right is smooth (we can identify $\Aut(\ms F(x))$ with $G$. This identification is not canonical, but the smooth structure on $\Aut(\ms F(x))$ is). We require that $\ms F$ preserves this smooth structure for some $x\in M$ or equivalently for all $x\in M$.
			
			\item Set $\Pi_2^\thin(M,x,y)$ to be the groupoid of paths $x\to y$ and bigons between them (with vertical composition of bigons as groupoid composition). Let $\Hom(\ms F(x)_0,\ms F(y)_1)$ be the set of equivariant morphisms $\ms F(x)_0\to \ms F(y)_1$ (in the sense that $f(p\cdot g) = f(p)\cdot \id_g)$. Then $\ms F$ induces a map $\Pi_2^\thin(M,x,y)\to \Hom(\ms F(x)_0,\ms F(y)_1)$. By picking $p\in \ms F(x)_0$ we identify $\Hom(\ms F(x)_0,\ms F(y)_1)\cong H$, $f\mapsto f(p):\id_{s(f(p))}$ (cf. lem. \ref{lem:g2tor2morphH}). Changing $p\in \ms F(x)_0$ changes this identification by conjugation, and therefore gives a well-defined smooth structure on $\Hom(\ms F(x)_0,\ms F(y)_1)$. Subsequently we require that $\Pi_2^\thin(M,x,y)\to \Hom(\ms F(x)_0,\ms F(y)_1)$ is smooth with respect to the diffeological structure on $\Pi_2^\thin(M,x,y)$ and the smooth structure on $\Hom(\ms F(x)_0,\ms F(y)_1)$ for some $x,y\in M$ (or equivalently for all $x,y\in M$).
		\end{itemize}

		Note defined this way, a transport 2-functor induces a transport functor $\tra^1\colon \Pi_1^\thin(M)\to G\tor$ in the sense of \cite{CollierLermanWolbert}. They showed that such functors are the same as principal $G$-bundles with connection. Here `the same' means an equivalence of the relevant categories. We want an analogous result for transport 2-functors and principal $\mc G$-2-bundles with 2-connection. To this end we first define the 2-category of transport 2-functors by adapting the definition from Schreiber \& Waldorf \cite{SW11}. They worked instead with 2-functors $\Pi_2^\thin(M)\to B \mc G$, which only captures trivial 2-bundles. 
		
		Next we describe the morphisms of the 2-category of 2-transport functors. These definitions can also be used to define a 2-functor 2-category $\mathbf{Fun}^2(\ms C,\ms D)$ for any pair of 2-categories. \cite{Power}
		
		\begin{definition}\label{def:2trans1morph}
			A 1-morphism $\rho\colon \ms F\to \ms F'$ of transport 2-functors is a \textit{pseudonatural transformation}. That is, an assignment $x\mapsto \rho(x)\colon \ms F(x)\to \ms F'(x)$ such that the following diagram commutes up to a 2-morphism $\rho(\gamma)$ for all $\gamma\colon x\to y$:
			\begin{equation}\label{eq:PseudonaturalTransform}
				\begin{tikzcd}[column sep=4em, row sep =4em]
					\ms F(x) \ar["\ms F(\gamma)"]{r} \ar["\rho(x)"']{d} & \ms F(y) \ar["\rho(y)"]{d} \ar["\rho(\gamma)" description,shorten >=0.3cm,shorten <=0.3cm, Rightarrow]{dl}\\
					\ms F'(x) \ar["\ms F'(\gamma)"']{r} & \ms F'(y)
				\end{tikzcd}
			\end{equation}
			We require that $\rho$ preserves composition. That is, for all $\gamma\colon x\to y$ and $\gamma'\colon y\to z$ we require:
			\begin{equation}\label{eq:1morphcomposition}
				\begin{tikzcd}[column sep=4em, row sep =4em]
					\ms F(x) \ar["\ms F(\gamma)"]{r} \ar["\rho(x)"']{d} & \ms F(y) \ar["\rho(y)" description]{d} \ar["\ms F(\gamma')"]{r} \ar["\rho(\gamma)" description,shorten >=0.3cm,shorten <=0.3cm, Rightarrow]{dl} & \ms F(z) \ar["\rho(z)"]{d} \ar["\rho(\gamma')" description,shorten >=0.3cm,shorten <=0.3cm, Rightarrow]{dl} \\
					\ms F'(x) \ar["\ms F'(\gamma)"']{r} & \ms F'(y) \ar["\ms F'(\gamma')"']{r} & \ms F'(z) 
				\end{tikzcd}
			    \quad=\quad
			    \begin{tikzcd}[column sep=4em, row sep =4em]
			    	\ms F(x) \ar["\ms F(\gamma'\circ \gamma)"]{r} \ar["\rho(x)"']{d} & \ms F(z) \ar["\rho(z)"]{d} \ar["\rho(\gamma'\circ\gamma)" description,shorten >=0.3cm,shorten <=0.3cm, Rightarrow]{dl}\\
			    	\ms F'(x) \ar["\ms F'(\gamma'\circ\gamma)"']{r} & \ms F'(z).
			    \end{tikzcd}
			\end{equation}
			We also require that $\rho$ is compatible with 2-morphisms. That is, for any $\Sigma\colon \gamma\Rightarrow \gamma'$ we require: 
			\begin{equation}\label{eq:1morph2morph}
				{\begin{tikzcd}[column sep=6em, row sep =4em]
					\ms F(x) \ar["\ms F(\gamma)"]{r} \ar["\rho(x)"']{d} & \ms F(y) \ar["\rho(y)"]{d} \ar["\rho(\gamma)" description,shorten >=0.3cm,shorten <=0.3cm, Rightarrow]{dl}\\
					\ms F'(x) \ar["\ms F'(\gamma)"{name=U,description}]{r} \ar[bend right=80, "\ms F'(\gamma')"{name=D,below}]{r} \ar[Rightarrow, shorten >=0.1cm, from=U, to=D,"\ms F'(\Sigma)"{description, pos=0.4}]& \ms F'(y) 
				\end{tikzcd}}
				\quad=\quad
				{\begin{tikzcd}[column sep=6em, row sep =4em]
					\ms F(x) \ar["\ms F(\gamma')"{name=D,description}]{r} \ar["\rho(x)"']{d} \ar[bend left = 80,"\ms F(\gamma)"{name=U}]{r} \ar[Rightarrow, from=U, to=D,shorten <=0.1cm, "\ms F(\Sigma)"{description,pos=0.5}]& \ms F(y) \ar["\rho(y)"]{d} \ar["\rho(\gamma')" description,shorten >=0.3cm,shorten <=0.3cm, Rightarrow]{dl}\\
					\ms F'(x) \ar["\ms F'(\gamma')"']{r} & \ms F'(y).
				\end{tikzcd}}
			\end{equation}
			Composition of two modifications $\rho:\ms F\to \ms F'$, $\rho':\ms F'\to \ms F''$ is defined by:
			\begin{align}
				(\rho'\circ \rho)(x) &= \rho'(x)\circ \rho(x),\\
				(\rho'\circ\rho)(\gamma) &= (\rho'(\gamma)\circ1_{\rho(x)})\circ (1_{\rho'(y)}\circ \rho(\gamma)),\qquad\text{for } \gamma:x\to y.
			\end{align}
		\end{definition}
	
		\begin{definition}\label{def:2trans2morph}
			A 2-morphism $\rho,\rho'$ of pseudonatural transformations of transport 2-functors $\ms F,\ms F'$ is a \textit{modification} $\mc A\colon \rho\Rightarrow\rho'$. That is, for each $x\in M$ a 2-morphism of $\mc G$-2-torsors
			\begin{equation}
				\begin{tikzcd}[column sep = 5em]
					\ms F(x) \ar["\rho(x)"{name=U},bend left=40]{r} \ar["\rho'(x)"{name=D,below}, bend right = 40]{r} \ar[Rightarrow, shorten >=0.1cm, shorten <=0.1cm, from=U, to = D, "\mc A(x)"{description}]& \ms F'(x)
				\end{tikzcd}
			\end{equation}
			Such that for each $\gamma\colon x\to y$ we have
			\begin{equation}\label{eq:2morphtransformrho}
				\begin{tikzcd}[column sep=4em, row sep =6em]
					\ms F(x) \ar["F(\gamma)"]{r} \ar["\rho'(x)"{left,name=L},bend right=100]{d} \ar["\rho(x)"{description,name=R}]{d} \ar[from =R, to=L, Rightarrow, "\mc A(x)"{above,pos=0.4}, shorten >= 0.1cm] & \ms F(y) \ar["\rho(y)"]{d} \ar["\rho(\gamma)" description,shorten >=0.3cm,shorten <=0.3cm, Rightarrow]{dl}\\
					\ms F'(x) \ar["\ms F'(\gamma)"']{r} & \ms F'(y)
				\end{tikzcd}
				\quad=\quad
				\begin{tikzcd}[column sep=4em, row sep =6em]
					\ms F(x) \ar["\ms F(\gamma)"]{r} \ar["\rho'(x)"']{d} & \ms F(y) \ar["\rho'(y)"{description, name = L}]{d} \ar["\rho(y)"{right, name=R}, bend left = 100]{d} \ar["\mc A(y)"{above,pos=0.5},from=R, to=L,Rightarrow,shorten <=0.1cm] \ar["\rho'(\gamma)" description,shorten >=0.3cm,shorten <=0.3cm, Rightarrow]{dl}\\
					\ms F'(x) \ar["\ms F'(\gamma)"']{r} & \ms F'(y)
				\end{tikzcd}
			\end{equation}
			Horizontal and vertical composition of 2-morphisms is given by pointwise horizontal and vertical composition of the corresponding 2-morphisms of $\mc G$-2-torsors.
		\end{definition}
		
		\begin{definition}
			2-transport functors $\Pi_2^\thin(M)\to\mc G\tor$ with pseudonatural transformations and modifications define a strict 2-category  $\Trans^2(M,\mc G)$.
		\end{definition}

\section{Principal $\mc G$-2-bundles with 2-connection}
The main objective of this paper is to give a useful definition of a principal $\mc G$-2-bundle with 2-connection based on the natural notion of 2-transport functor. It is not clear a priori what such a 2-bundle with 2-connection is, and it is also not obvious what the right notion of morphism are for such bundles. We will give a definition of 2-bundles with 2-connection that is constructed in such a way that it will give a 2-category equivalent to the 2-category of 2-transport functors. We will skip ahead and give the definition, and then in the remainder of the section we justify this definition by showing its equivalence to 2-transport functors. 

	\begin{definition}
		A principal $\mc G$-2-bundle  over a manifold $M$ is a groupoid $\mc P_1\rightrightarrows \mc P_0$  with a right action $R\colon \mc P\circlearrowleft \mc G$ and a $\mc G$-invariant surjective submersion functor $\pi\colon \mc P\to (M\rightrightarrows M)$, i.e. surjective submersions fitting in a commutative diagram:
		\[
		\begin{tikzcd}[column sep = 1ex]
		{G{\ltimes} H}\circlearrowright\mc P_1 \ar[rr,shift left=0.5ex]\ar[rr,shift right=0.5ex]\ar[dr,"\pi_1"',start anchor={[xshift=2.2ex]}] &&\mc P_0 \circlearrowleft G\ar[dl,"\pi_0",start anchor = {[xshift=-1ex]}]\\
		&M
		\end{tikzcd}
		\]
		Moreover we require that the following functor is an isomorphism of categories (i.e. there is an inverse)
		\begin{equation}\label{eq:princBundleAxiom}
			(\mr{pr}_1,R)\colon \mc P\times\mc G\to \mc P\times_M\mc P,\qquad (p,g)\mapsto (p,p\cdot g),
		\end{equation}
		where $\mr{pr}_1$ is projection to the first factor.
	\end{definition}
	
	\begin{remark}
	Since $\mc P$ is a groupoid, there is an identity map $\mc P_0\to \mc P_1$, therefore a trivialization of $\mc P_0$ induces a trivialization of $\mc P_1$ (we can extend any local section of $\mc P_0$ to a section of $\mc P_1$, alternatively transition functions of $\mc P_0$ define transition functions of $\mc P_1$ under the inclusion $G\to G{\ltimes}H$). Furthermore the definition of a 2-connection and the definitions of all the relevant morphisms can be completely phrased in terms of $\mc P_0$. Therefore in this setting we can completely replace principal $\mc G$-2-bundles by ordinary principal $G$-bundles. However in the proofs the groupoid structure and $\mc G$ action naturally appear, and therefore it is better to keep this definition as it is.
	\end{remark}

	\begin{definition}
		A 2-connection on a principal $\mc G$-2-bundle $\mc P_1\rightrightarrows \mc P_0$ is a pair $(A,B)$ with $A\in \Omega^1(\mc P_0,\mf g)$ a connection on $\mc P_0$ and $B\in \Omega^2(\mc P_0,\mf h)$ satisfying:
		\begin{enumerate}
		\item (Equivariance): $R_g^*B = (\alpha_g^{-1})_*B$ for all equivariant $g\colon \mc P_0\to G$. Here $(\alpha_g^{-1})_*$ is obtained by differentiating $\alpha_{g(m)}^{-1}\colon H\to H$ to a map $\mf h\to \mf h$ for each $m\in M$.
		\item (Fake-flatness): $t_*B=F_A$ where $t_*\colon \mf h\to \mf g$ and $F_A$ is the curvature of $A$.
		\end{enumerate}
The equivariance condition is equivalent to 
	\begin{enumerate}
	 \item (Equivariance): $R_g^*B = (\alpha_g^{-1})_*B$ for all $g\in G$ and $\iota(X_\xi)B=0$ for all $\xi \in \mf g$, where $X_\xi$ is the fundamental vector field of the $G$ action on $\mc P_0$.	
	\end{enumerate}
	\end{definition}

This definition of principal 2-bundle with 2-connection coincides with that of 'special $\mc G$ 2-bundles' in \cite{MartinsPicken1}, however they do not explain the equivalence of that definition with 2-transport. In the category of principal bundles with connections one requires that morphisms preserve the connection. That is, a 1-morphism $\rho\colon (P,A)\to (P',A')$ is an equivariant bundle map such that $\rho^*A'=A$. On the side of transport functors this is the requirement that a 1-morphism is a natural transformation and hence commutes with parallel transport in the sense of diagram \eqref{eq:PseudonaturalTransform} (on the nose). However for transport 2-bundles we require this diagram only to commute up to 2-morphism. This translates to a slightly different notion of 1-morphisms for 2-bundles with 2-connection.
	
	\begin{definition}
		A 1-morphism of 2-bundles with 2-connection $(\mc P,A,B)\to (\mc P',A',B')$ is a pair $(F,\phi)$ with $F\colon \mc P\to \mc P'$ an equivariant functor (or equivalently an equivariant bundle map $\mc P_0\to \mc P'_0$), such that $\pi'\circ F=\pi$ and $\phi\in \Omega^1(\mc P_0,\mf h)$ a form such that:
		\begin{enumerate}
			\item $R_g^*\phi=(\alpha_g^{-1})_*\phi$ for all equivariant $g\colon \mc P_0\to G$, or equivalently, \\$R_g^*\phi=(\alpha_g^{-1})_*\phi$ for all $g\in G$ and $\phi(X_\xi)=0$ for all $\xi\in \mf g$.
			\item $F^*A' = A+t_*\phi$.
			\item $F^*B' = B+d\phi+\frac12[\phi,\phi]+\frac12[A\wedge \phi]$, where $[A\wedge\phi]$ is understood as the Lie bracket of $\mf g\ltimes \mf h$-valued forms (cf. eq. \eqref{eq:infcrossedmodulebracket}). With slight abuse of notation one can replace $\frac12[A\wedge \phi]$ by $\alpha_*(A\wedge\phi)$.
		\end{enumerate}
		Composition is defined by $(F',\phi')\circ (F,\phi) = (F'\circ F,\phi+F^*\phi')$.  
	\end{definition}
	
 	We can also consider automorphisms, and take $R_g\colon \mc P_0\to \mc P_0$ for equivariant $g\colon \mc P_0\to G$ (and still $\phi\in \Omega^1(\mc P_0,\mf h)$). Then we get the appropriate notion of gauge transformations for 2-connections:
	\begin{enumerate}
	\item $A\mapsto \Ad_g^{-1}\left(A+dgg^{-1}+t_*\phi\right)$ 
	\item $B\mapsto (\alpha_g^{-1})_*	\left(B+d\phi+\frac12[\phi,\phi]+\frac12[A\wedge \phi]\right)$
	\end{enumerate}
	In physical terms, the action of $\phi$ should be interpreted as a higher gauge symmetry.

	\begin{definition}
		A 2-morphism of 2-bundles with 2-connection 
		\begin{equation}
			\begin{tikzcd}[column sep=2em]
				(\mc P, A, B) \ar["{(F,\phi)}"{name=U},bend left= 50]{r} \ar[bend right = 50, "{(F',\phi')}"{name = D,below}]{r} \ar[from= U,to=D,Rightarrow,"a"{description},shorten <=1mm,shorten >=1mm] & (\mc P',A',B')
			\end{tikzcd}
		\end{equation}
		is a map $a\colon \mc P_0\to H$ such that for all $p\in \mc P_0$:
		\begin{align}
			a(p\cdot g)&=\alpha_{g^{-1}}a(p)\\
			F'(p)&=t(a(p))F(p)\\
			\phi' & = \Ad_{a}\phi-(r_a^{-1}\circ \alpha_a)_*A-a^*\theta
		\end{align}
		Here $r_{a}^{-1}\circ \alpha_{a}$ is the map $\mc P_0\times G\to H$, $(p,g)\mapsto \alpha(g,a(p))a(p)^{-1}$. Vertical composition of $a,a'$ is pointwise multiplication $a(p)\cdot a'(p)$. Horizontal composition of $a\colon (F,\phi)\Rightarrow(F',\phi')$, $a'\colon (\tilde F,\tilde \phi)\Rightarrow (\tilde F',\tilde \phi')$ is given by $(F')^*a'\cdot a$.
	\end{definition}
	
	This defines a strict 2-category of principal $\mc G$-2-bundles with 2-connection. One needs to check that this data satisfies all the axioms of a 2-category. This is a straightforward task, and we present it without proof.

	\begin{proposition}
		Principal $\mc G$-2-bundles with 2-connections and the 1- and 2-morphisms described above define a 2-category $\Bun_\nabla^2(M,\mc G)$.
	\end{proposition}
	
	Let $\mc P_1 \rightrightarrows \mc P_0$ be a principal $\mc G$-2-bundle. Then a trivialization of $\mc P_0$ with respect to some cover $U_i$ induces a trivialization of $\mc P_1$ (both as ordinary principal bundles) by using the identity map $\mc P_0\to \mc P_1$ to extend any section of $\mc P_0$. Alternatively if we have transition functions $g_{ij}\colon U_{ij}\to G$ then we can trivially extend these to transition functions $U_{ij}\to G\ltimes H$. With respect to such a trivialization a 2-connection is given by
	\begin{equation*}
		A_i\in \Omega^1(U_i,\mf g),\qquad B_i\in \Omega^2(U_i,\mf h),
	\end{equation*}
	such that:
	\begin{align}
		t_*B_i &= F_i = dA_i+\frac12[A_i,A_i],\\
		A_i &= g_{ij}^{-1}A_jg_{ij} + g_{ij}^{-1}dg_{ij},\\
		B_i &= (\alpha_{g_{ij}}^{-1})_*B_j.
	\end{align}

\section{Equivalence between transport 2-functors and principal {2-bundles} with 2-connection}

	\subsection{Extracting geometric data from transport 2-functors}
	
	This section will be devoted to showing the equivalence between transport 2-functors and principal 2-bundles with 2-connection. To start, we will show how to extract the data of a principal 2-bundle with 2-connection from a transport 2-functor. 
	
	First we will recall how to obtain a principal $G$-bundle from a smooth functor $\tra:\Pi_1^\thin(M)\to \mc G\tor$, following \cite{CaetanoPicken}. Let us fix some $*\in M$ then let $\Pi_1^\thin(M,*)\subset \Pi_1^\thin(M)$ be the full subgroupoid of paths starting at $*\in M$, and let $\Omega_1^\thin(M,*)\subset \Pi_1^\thin(M,*)$ be the group of thin homotopy classes of loops at $*$. Then fixing some isomorphism $\tra(*)\cong G$ we obtain a homomorphism $\hol:\Omega_1^\thin(M,*)\to G$. We then consider the following equivalence relation on $\Pi_1^\thin(M,*)\times G$:
	\begin{equation}
	 (\gamma,g)\sim (\gamma',h) \quad \text{if} \quad \gamma(1)=\gamma'(1),\quad \hol(\gamma'\gamma^{-1})=hg^{-1}.
	\end{equation}
	Then consider $P=\Pi_1^\thin(M,*)\times G/\sim$. With projection $\pi[\gamma,g]=\gamma(1)$ and $G$-action $[\gamma,g] \cdot h= [\gamma,gh]$ this defines a principal $G$-bundle $P\to M$. 
	
	Recall that parallel transport along a path $\gamma\colon x\to y$ is defined as follows. Let $\tilde \gamma$ be (any) lift of $\gamma$, starting at some $p\in P_x$. Then we have
	\begin{equation}
		\tra_\gamma(p) : \tilde\gamma(1) = \mc P\exp \int_0^1\!-A(\tilde \gamma'(t))\,\dd t,
	\end{equation}
	where the right hand side is the path ordered exponent, i.e. the solution to the differential equation
	\begin{equation}
		g'(t) = -A(\gamma'(t))g(t),\qquad g\colon I\to G.
	\end{equation}
	This allows us to recover the connection from parallel transport. Let $\gamma_t$ for $t\in I$ be the path $\gamma_t(s)=\gamma(st)$, then note that
	\begin{equation}\label{eq:Adef}
		A(\tilde\gamma'(0)) = -\left.\frac{\dd }{\dd t}\right|_{t=0}\tra(\gamma_t)(p):\tilde \gamma(t).
	\end{equation}
	Or alternatively,
	\begin{equation}
		A(\tilde\gamma'(0)) = \left.\frac{\dd }{\dd t}\right|_{t=0}\tra(\gamma_t)^{-1}(\tilde\gamma(t)):\tilde\gamma(0).
	\end{equation}
	The smoothness condition on transport functors then assures that this defines a $\mf g$-valued differential form. 
	
	On the other hand given a principal $G$-bundle with connection, we obtain a smooth functor $\Pi_1^\thin(M)\to \mc G\tor$ through parallel transport. These two constructions give an equivalence of categories. Stated in different terms this result is originally due to \cite{Barrett}, but using categorical language this statement appears in \cite{SW09, CollierLermanWolbert}.
	\begin{theorem} \label{thm:transportisbundle}
		The procedure above defines an equivalence of categories $\Bun_\nabla^1(M,G)\cong \Trans^1(M,G)$, i.e. respectively the category of principal $G$-bundles with connection on $M$ and the category of smooth functors $\Pi_1^\thin(M)\to G\tor$.
	\end{theorem}
	
\begin{figure}[!htb]\centering
			\includegraphics[width=0.3\textwidth]{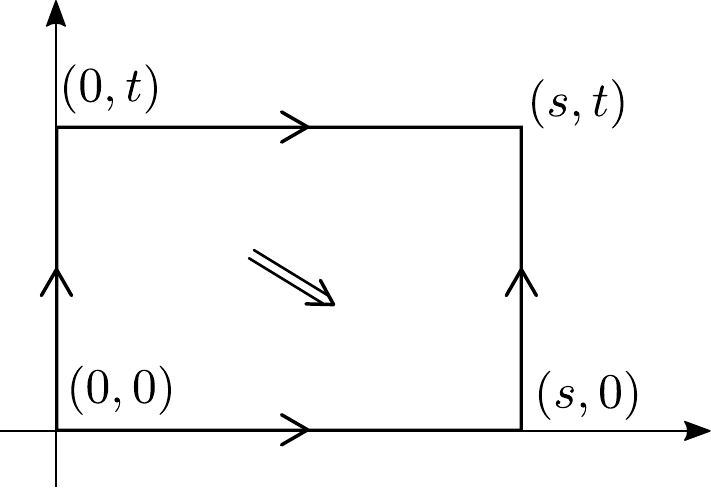}
			\caption{\label{fig:canonicalBigon}Canonical bigon $\Gamma(s,t)$ in $I^2$. Any two bigons in $I^2$ with this source and target are necessarily thinly homotopic since by dimensional reasons any homotopy of bigons is thin. A simple parameterisation is given by $\Gamma(s,t)(u,v)=(sv,tv^{\frac{1-u}{u}})$.
			}
		\end{figure}

	\noindent The $B$-form is obtained from 2-transport in a similar, but more involved way. Let $\Sigma\colon I^2\to P$ be a smooth map, and denote
		\begin{equation}
			X=\left.\frac{\dd}{\dd s}\right|_{s=0}\Sigma(s,0),\qquad Y=\left.\frac{\dd}{\dd t}\right|_{t=0}\Sigma(0,t).
		\end{equation}
		To each pair $(s,t)\in I^2$ we can associate a canonical thin homotopy class of bigons $\Gamma(s,t):I^2\to I^2$. This is shown in figure \ref{fig:canonicalBigon}. This gives a map $\Gamma\colon I^2\to \Pi^\thin_2(I^2)^2$, i.e. the space of bigons on $I^2$. Composition with $\Sigma$ then defines a map $\Sigma_*\colon \Pi^\thin_2(I^2)\to \Pi^\thin_2(P)$. Thus $\Sigma_*\Gamma(s,t)$ is a bigon for each $s,t\in I^2$. We obtain a 2-form by taking the 2-transport along $\Sigma_*\Gamma(s,t)$ and differentiating with respect to $s,t$. More precisely, we define
		\begin{equation}\label{eq:Bdefinition}
			B(X,Y) = \left.\frac{\partial^2}{\partial s\partial t} \right|_{(0,0)} \tra^2(\pi_*\Sigma_*\Gamma(s,t))^{-1}( \Sigma(s,t)):\id_{\tra(s(\pi_*\Sigma_*\Gamma(s,t))^{-1}(\Sigma(s,t))}.
		\end{equation}
		Here $s$ is the source map, mapping bigons to their source path. This definition seems complicated at first glance. The point is that $\tra^2(\pi_*\Sigma_*\Gamma(s,t))^{-1}( \Sigma(s,t))$ always lies in the same fiber as $\id_{\Sigma(0,0)}$. Therefore we could in principle divide by $\id_{\Sigma(0,0)}$ instead. However if we did that, we would end up with a $\mf g\ltimes \mf h$-valued form. Instead we divide by parallel transport of the source of $\pi_*\Sigma_*\Gamma(s,t)$ to get something $\mf h$-valued. This is precisely the same as dividing by $\id_{\Sigma(0,0)}$ and then applying the projection $\mf g\ltimes\mf h\to \mf h$. 
		
		\begin{theorem}
			The definition \eqref{eq:Bdefinition} of $B$ does not depend on choice of $\Sigma\colon I^2\to P$ (fixing $X,Y\in TP$), and defines a differential form $B\in \Omega^2(P,\mf h)$ satisfying fake flatness $t_*B=F_A$ and the following equivariance condition:
			\begin{equation}
				R_g^*B = (\alpha_g^{-1})_*B,
			\end{equation}
			for any equivariant $g\colon \mc P_0\to G$, where $(\alpha_g)_*\colon \mf h\to \mf h$ is obtained by differentiating the $G$ action on $H$.
		\end{theorem}	
		\proof This theorem is proved for trivial 2-bundles by Schreiber and Waldorf \cite{SW11}, see also \cite{Parzygnat}. This proof is rather technical, but does provide some insight. Since 2-bundles admit local trivializations and the statement of the theorem is local, it only remains to show the equivariance condition. This is a straightforward computation:
		\mathindent2mm
		\begin{align*}
			R_g^*B(X,Y) &= \left.\frac{\partial^2}{\partial s\partial t} \right|_{(0,0)} \tra^2(\pi_*\Sigma_*\Gamma(s,t))^{-1}(\Sigma(s,t)\cdot g(\Sigma(s,t)):\id_{\tra(s(\pi_*\Sigma_*\Gamma(s,t))^{-1}(\Sigma(s,t))\cdot g(\Sigma(s,t))}\\
			&=\left.\frac{\partial^2}{\partial s\partial t} \right|_{(0,0)} \tra^2(\pi_*\Sigma_*\Gamma(s,t))^{-1}(\Sigma(s,t))\cdot 1_{g(\Sigma(s,t))}:\left[\id_{\tra(s(\pi_*\Sigma_*\Gamma(s,t))^{-1}(\Sigma(s,t))}\cdot 1_{g(\Sigma(s,t))}\right]\\
					&=\left.\frac{\partial^2}{\partial s\partial t} \right|_{(0,0)} 1_{g(\Sigma(s,t))^{-1}}\cdot \left[\tra^2(\pi_*\Sigma_*\Gamma(s,t))^{-1}(\Sigma(s,t)):\id_{\tra(s(\pi_*\Sigma_*\Gamma(s,t))^{-1}(\Sigma(s,t))}\right]\cdot 1_{g(\Sigma(s,t))}\\
			&=\left.\frac{\partial^2}{\partial s\partial t} \right|_{(0,0)} \alpha_g^{-1}\left[\tra^2(\pi_*\Sigma_*\Gamma(s,t))^{-1}(\Sigma(s,t)):\id_{\tra(s(\pi_*\Sigma_*\Gamma(s,t))^{-1}(\Sigma(s,t))}\right]\\
			&=(\alpha_g^{-1})_*B(X,Y).\qedhere
		\end{align*}\mathindent12.5mm
		
		\subsection{The parallel 2-transport of a 2-connection}
		Given a transport 2-functor, we now know how to obtain a principal bundle with 2-connection $(A,B)$. The main tool to understand how to obtain transport 2-functors from a 2-connection is the non-Abelian Stokes Theorem. This theorem expresses the transport around a contractible loop as an integral of the curvature over a disk bounding the loop. 
		
		\begin{theorem}[\textbf{non-Abelian Stokes}]\label{thm:nonabelianstokes}
			Let $A$ be a connection on a principal $G$-bundle $P$, and let $\Gamma_t$ be a bigon. Let $F_A$ be the curvature of $A$, and let $(\partial s,\partial t)$ be the natural global frame on $I^2$. Then
			\begin{equation}
				\tra(\Gamma_1)(p):\tra(\Gamma_0)(p) = \mc P\exp \int_0^1\!\!\dd s\left[\int_0^1\!\! \dd t \,\,\tilde{\Gamma}^*F_A(\partial s,\partial t)\right],
			\end{equation}
			where $\tilde{\Gamma}\colon I^2\to P$ is obtained by horizontally lifting $\Gamma_s$ at $p$ for each $s$ individually, i.e. 
			\[
				\tilde\Gamma(s,t) = \tra(\Gamma_{s,t})(p).
			\]
		\end{theorem}
		
\begin{figure}[!htb]\centering
			\includegraphics[width = 0.9\textwidth]{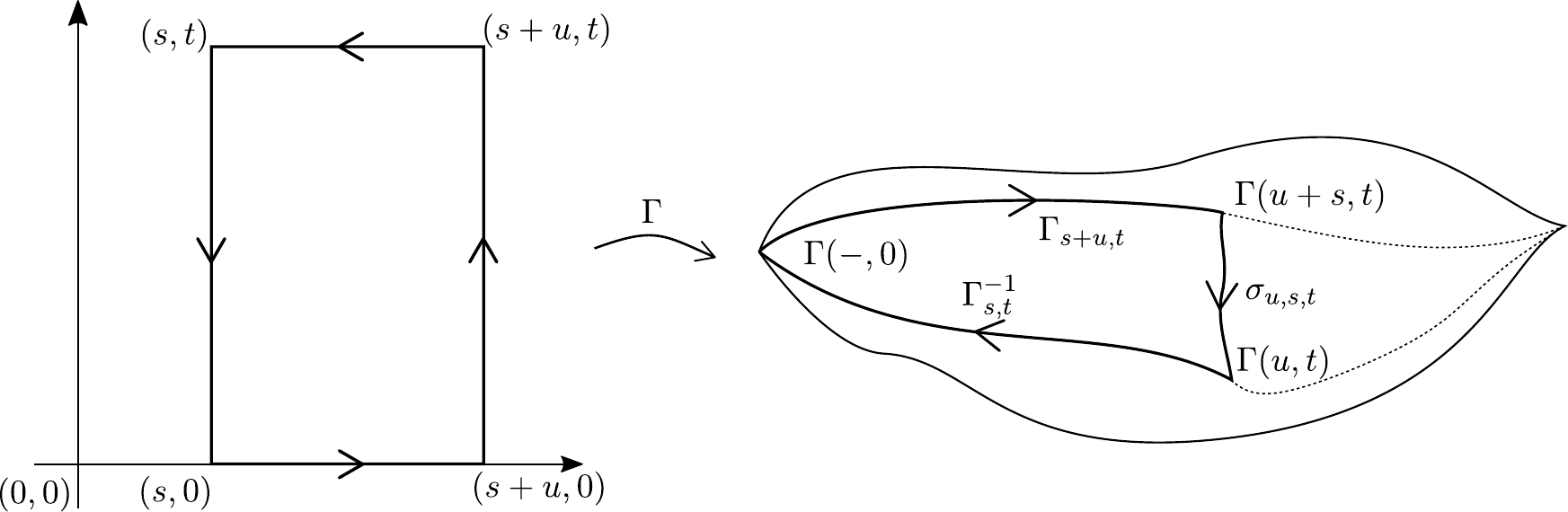}
			\caption{Sketch of the loop used in the proof of the non-Abelian Stokes Theorem.\label{fig:nonAbelianStokesLoop}}
		\end{figure}

		\proof Let $\Gamma_u$ denote the path $\Gamma(u,\cdot)$, and let $p\in P_x=P_{\Gamma(\cdot,0)}$ be fixed. Then consider
		\[
			f(s) = \tra({\Gamma_0})(p):\tra(\Gamma_s)(p).
		\]
		This function satisfies
		\[
			f(s+u)f(s)^{-1}=\tra(\Gamma_{s})(p):\tra(\Gamma_{s+u})(p).
		\]
		Thus differentiating we get the following differential equation:
		\begin{equation}
			f'(s) = -\left[\left.\frac{\partial}{\partial u}\right|_{u=0}\tra(\Gamma_{s+u})(p):\tra(\Gamma_s)(p)\right]\cdot f(s).
		\end{equation}
		The solution of this differential equation is a path ordered exponential. Now we will introduce an additional parameter. Let $\Gamma_{s,t}$ denote the path $\tau\mapsto \Gamma(s,t\tau)$, and let \[\sigma_{u,s,t}(\tau) = \Gamma(u+s\tau,t).\]
		 Finally denote $\tilde \Gamma(s,t) = \tra(\Gamma_{s,t})(p)$. Then consider the parallel transport along the following loop (see also figure \ref{fig:nonAbelianStokesLoop})
		\begin{align*}
			\tra(\Gamma_{s,t}&)^{-1}\circ \tra(\sigma_{u,s,t})^{-1}\circ \tra(\Gamma_{u+s,t})(p):p\\
			&=\tra(\sigma_{u,s,t})^{-1}\left(\tilde\Gamma(s+u,t)\right):\tilde\Gamma(s,t).
		\end{align*}
		
		\noindent For $t=0$ this loop is constant, and for $t=1$ this is $f(s)^{-1}f(s+u)$. Therefore we deduce the equality
		\begin{equation*}
			\tra(\Gamma_0)(p):\tra(\Gamma_1)(p) = \mc P\exp \int_0^1\!\!\dd s\left[\int_0^1\!\! \dd t -\left.\frac{\partial^2}{\partial u\partial t}\right|_{(0,t)} \left(\tra(\sigma_{u,s,t})^{-1}\left(\tilde\Gamma(s+u,t)\right):\tilde\Gamma(s,t)\right) \right].
		\end{equation*}
		To complete the proof we just need to equate the integrand with the curvature. We note that
		\begin{align*}
			\left.\frac{\partial^2}{\partial u\partial t}\right|_{(0,t)}& \left(\tra(\sigma_{u,s,t})^{-1}\left(\tilde\Gamma(s+u,t)\right):\tilde\Gamma(s,t)\right)\\
			 &= \frac{\partial}{\partial t}A\left(\left.\frac{\partial}{\partial u}\right|_0\tilde \Gamma(s+u,t)\right)\\
			 &= \left(\mc L_{\tilde \Gamma_*\partial t} A\right)\left(\tilde\Gamma_*\partial s\right)\\
			\intertext{Using Cartan's magic formula and $A(\tilde\Gamma_*\partial t)=0$ we obtain}
			&=dA\left(\tilde \Gamma_*\partial s,\tilde\Gamma_*\partial t\right).
		\end{align*}
		For a vector field $X$ let $X^h$ and $X^v$ denote the horizontal and vertical part of $X$ respectively, with respect to the connection $A$. Then recall that the curvature of $A$ is given by
		\[
		 	F(X,Y)=dA(X^h,Y^h).
		\]
		Since $\Gamma_*\partial t$ is horizontal, we are done if we show that
		\[
			dA(\tilde \Gamma_*\partial s^v,\tilde\Gamma_*\partial t) = 0.
		\]
		This is a simple computation,
		\[
			dA(\tilde \Gamma_*\partial s^v,\tilde\Gamma_*\partial t) = \tilde \Gamma_*\partial s^vA(\tilde \Gamma_*\partial t)-\tilde \Gamma_*\partial tA(\tilde \Gamma_*\partial s^v)-A([\tilde \Gamma_*\partial s^v,\tilde \Gamma_*\partial t]).
		\]
		The first term is zero because $\tilde \Gamma_*\partial t$ is horizontal. The second term is zero because $\tilde \Gamma_*\partial t$ is horizontal and $A(\tilde \Gamma_*\partial s^v)$ is vertical and hence live in complementary subspaces of $TP$. Finally the last term is zero because the flows of $\tilde \Gamma_*\partial t$ and $\tilde \Gamma_*\partial s^v$ commute (and hence their Lie bracket vanishes). \qed
		
		We will now use the non-Abelian Stokes Theorem to derive a formula for 2-transport. The 2-transport should assign to a bigon $\Gamma_s\colon x\to y$ a 2-morphism of $\mc G$-torsors $\tra(\Gamma_0)\Rightarrow\tra(\Gamma_1)$. Recall by lemma \ref{lem:g2tor2morphH} that for $p\in P_x$ we should therefore obtain an element  of $H$ such that
			\begin{equation}
				t\left(\tra^2_H(p):\id_{\tra(\Gamma_0)(p)}\right) = \tra^1(\Gamma_1)(p):\tra^1(\Gamma_0)(p).
			\end{equation}
			The right hand side of this equation appears in the non-Abelian Stokes Theorem. Now note that the right hand side of the non-Abelian Stokes Theorem has a natural preimage in $H$; we can simply replace $F$ by $B$.
			Therefore we can define 2-transport by:
			\begin{equation}\label{eq:2transportFormula}
				\tra^2(\Gamma)(p) = \id_{\tra^1(\Gamma_0)(p)}\cdot \,\,\mc P\exp \int_0^1\!\!\dd s\left[\int_0^1\!\! \dd t \left(\tilde{\Gamma}_s^*B(\partial s,\partial t)\right)\right].
			\end{equation}
			This definition satisfies the equivariance condition of lemma \ref{lem:g2tor2morphH} as well:
			\begin{align*}
				\tra^2(\Gamma)(p\cdot g) &= \id_{\tra^1(\Gamma_0)(p\cdot g)}\cdot \,\,\mc P\exp \int_0^1\!\!\dd s\left[\int_0^1\!\! \dd t \left(\tilde{\Gamma}_s^*R_g^*B(\partial s,\partial t)\right)\right]\\
				&= \id_{\tra^1(\Gamma_0)(p\cdot g)}\cdot \,\,\mc P\exp \int_0^1\!\!\dd s\left[\int_0^1\!\! \dd t\, (\alpha_g^{-1})_*\left(\tilde{\Gamma}_s^*B(\partial s,\partial t)\right)\right]\\
				&= \id_{\tra^1(\Gamma_0)(p\cdot g)}\cdot \,\,\alpha_g^{-1}\left(\mc P\exp \int_0^1\!\!\dd s\left[\int_0^1\!\! \dd t \left(\tilde{\Gamma}_s^*B(\partial s,\partial t)\right)\right]\right).
			\end{align*}
			
			This expressions defines a transport 2-functor (cf. def. \ref{def:transport2functor}). This is shown in \cite{SW11} for trivial bundles, and extended in \cite{Voorhaar} for non-trivial bundles. The compatibility with vertical and horizontal composition is a relatively simple computation. The most involved part is proving thin-homotopy invariance. Thin homotopy invariance follows from the higher non-Abelian Stokes theorem \ref{thm:nonabelianstokes2}, however the proof of this theorem uses the vertical and horizontal composition rules for bigons. Therefore we could fix some particular way of parameterizing horizontal and vertical composition of bigons, then prove the higher non-Abelian Stokes theorem and use it to conclude that the result is independent of this parameterization. In summary, we have:
			
			\begin{theorem}
				Equation \eqref{eq:2transportFormula} defines a 2-transport functor.
			\end{theorem}
			
		\subsection{Equivalence on the level of morphisms}
		Now we know how to obtain a principal 2-bundle with 2-connection from a transport 2-functor and visa versa. We want to upgrade this to an equivalence of 2-categories. We do this by constructing two 2-functors. We will put more emphasis on the 2-functor taking transport 2-functors to principal 2-bundles with 2-connection, since the main objective of this paper is to justify our definition of principal 2-bundles with 2-connection.

		Let $\ms F,\ms G$ be two transport 2-functors, respectively defining a principal $\mc G$-2 bundle with 2-connection $(\mc P,A,B),\, (\mc P',A',B')$. Given a 1-morphism $\rho\colon\ms F\to \ms G$, we in particular obtain for each $x\in M$ a 1-morphism of $\mc G$-torsors $\rho(x)\colon \ms F(x)\to \ms G(x)$. This defines a bundle map $\colon\mc P\to \mc P'$. Then furthermore we obtain a 2-morphism $\rho(\gamma)$ of $\mc G$-torsors for each path $\gamma\colon x\to y$ in $M$, fitting in a diagram
		\begin{equation}
			\begin{tikzcd}[column sep=4em, row sep =4em]
				\mc P_{x} \ar["\tra^A_\gamma"]{r} \ar["F"']{d} & \mc P_y \ar["F"]{d} \ar["\rho(\gamma)" description,shorten >=0.3cm,shorten <=0.3cm, Rightarrow]{dl}\\
				\mc P'_x \ar["\tra^{A'}_\gamma"']{r} & \mc P'_y
			\end{tikzcd}
		\end{equation}
		Now giving a path $\gamma$ in $P$, we define:
		\begin{equation}\label{eq:phidef}
			\phi(\gamma'(0))=-\left.\frac{\dd}{\dd t}\right|_{t=0}\rho(\pi_*\gamma_t)(\gamma(0)):F\left(\tra_\gamma^A(\gamma(0))\right).
		\end{equation}
		Compare this to equation \eqref{eq:Adef} defining the connection. We claim the pair $(F,\phi)$ defines a 1-morphism $(\mc P,A,B)\to (\mc P',A',B')$. Conversely given a 1-morphism $(F,\phi)\colon(\mc P,A,B)\to (\mc P',A',B')$ we obtain a 1-morphism of $\mc G$-torsors $F(x):\mc P_x\to \mc P_x'$, and for every path $\gamma\colon x\to y$ in $M$ and $p\in \mc P_{x,0}$ we obtain
		\begin{equation}\label{eq:rhodef}
			\rho(\gamma)(p):\id_{F(\tra_\gamma(p))} = \mc P\exp \int_0^1\!\dd t\,\phi(\tilde \gamma'(t)),
		\end{equation}
		where $\tilde \gamma$ is the horizontal lift of $\gamma$ starting at $p$, and we claim $F,\rho$ defines a 1-morphism of the transport 2-functors associated to $(\mc P,A,B)$ and $(\mc P',A',B')$. 
	
		\begin{lemma}\label{lem:2functor1morphism}
			The pair $(F,\phi)$ where $\phi$ is defined by \eqref{eq:phidef} defines a 1-morphism of principal 2-bundles with 2-connection, i.e. $\phi\in \Omega^2(\mc P_0,\mf h)$ satisfying $R_g^*\phi=(\alpha_g^{-1})_*$ for all equivariant $g:\mc P_0\to G$ and 
			\begin{align}
			F^*A'&=A+t_*\phi,\\
			F^*B' &= B+d\phi+\frac 12[\phi\wedge\phi]+[A\wedge \phi].
			\end{align}
			Conversely $\rho$ defined by \eqref{eq:rhodef} defines a 1-morphism of transport 2-functors, i.e. it satisfies relations \eqref{eq:1morphcomposition} and \eqref{eq:1morph2morph}.
		\end{lemma}
	
		\proof  Let $\rho$ be a 1-morphism of transport 2-functors. In particular we obtain equivariant maps $\rho(x)\colon \ms F(x)\to \ms G(x)$ for all $x\in M$. By smoothness and lemma \ref{lem:equivfunctor} this is precisely a bundle map $F\colon \mc P_0\to \mc P'_0$. The further conditions (eq. \eqref{eq:PseudonaturalTransform}--\eqref{eq:1morph2morph}) give compatibility with transport.
		
		For any $\gamma\colon x\to y$ we have a 2-morphism:
		\begin{equation}\label{eq:defrhogamma}
			\begin{tikzcd}[column sep=4em, row sep =4em]
				\mc P_{x} \ar["\tra^A_\gamma"]{r} \ar["F"']{d} & \mc P_y \ar["F"]{d} \ar["\rho(\gamma)" description,shorten >=0.3cm,shorten <=0.3cm, Rightarrow]{dl}\\
				\mc P'_x \ar["\tra^{A'}_\gamma"']{r} & \mc P'_y
			\end{tikzcd}
		\end{equation}
		Using lemma \ref{lem:g2tor2morphH} we thus get a map $\rho_H(\gamma)\colon \mc (P_0)_x\to H$ for each $\gamma$. Recall that a 1-form $\phi\in \Omega^1(\mc P_0,\mf h)$ is the same as a transport functor $\Pi_1^\thin(\mc P_0)\to H$. Thus if $\rho_H$ is functorial with respect to composition of paths, we can differentiate to extract a 1-form $\phi$. One can check that $\rho_H$ is not functorial, but we can modify it to be functorial. To understand how to do this, let us consider the compatibility axiom of $\rho$ with path composition given by equation \eqref{eq:1morphcomposition}. Let $\gamma\colon x\to y$ and $\gamma'\colon y\to z$, then we have:
		\begin{equation}
			\begin{tikzcd}[column sep=4em, row sep =4em]
				\mc P_x \ar["\tra^A_\gamma"]{r} \ar["F"']{d} & \mc P_y \ar["F" description]{d} \ar["\tra^A_{\gamma'}"]{r} \ar["\rho(\gamma)" description,shorten >=0.3cm,shorten <=0.3cm, Rightarrow]{dl} & \mc P_z \ar["F"]{d} \ar["\rho(\gamma')" description,shorten >=0.3cm,shorten <=0.3cm, Rightarrow]{dl} \\
				\mc P'_x \ar["\tra^{A'}_\gamma"']{r} & \mc P'_y \ar["\tra^{A'}_{\gamma'}"']{r} & \mc P'_z 
			\end{tikzcd}
			\quad=\quad
			\begin{tikzcd}[column sep=4em, row sep =4em]
				\mc P_{x} \ar["\tra^A_{\gamma'\gamma}"]{r} \ar["F"']{d} & \mc P_z \ar["F"]{d} \ar["\rho(\gamma'\gamma)" description,shorten >=0.3cm,shorten <=0.3cm, Rightarrow]{dl}\\
				\mc P'_x \ar["\tra^{A'}_{\gamma'\gamma}"']{r} & \mc P'_z
			\end{tikzcd}
		\end{equation}
		Let $p\in \mc P_{x,0}$ then this means
		\begin{align}
			\rho(\gamma'\gamma)(p) &= \tra_{\gamma'}(\rho(\gamma)(p)) \circ \rho(\gamma')(\tra_\gamma(p)) \\
			&= \left(1_{\tra_{\gamma'}\circ F\circ \tra_{\gamma}(p)}\cdot\rho(\gamma)(p):1_{F\circ\tra_{\gamma}(p)}\right) \circ \rho(\gamma')(\tra_\gamma(p))
		\end{align}
		Here we used lemma \ref{lem:equivfunctor} to compute the first factor on the right hand side. We can divide by the source on both sides and use equation \eqref{eq:DivisionFunctorial} to get
		\begin{align}
			\rho(\gamma'\gamma)(p):1_{F\circ \tra_{\gamma'\gamma}(p)}&=\left(1_{\tra_{\gamma'}\circ F\circ \tra_{\gamma}(p)}\cdot\rho(\gamma)(p):1_{F\circ\tra_{\gamma}(p)}\right):1_{F\circ \tra_{\gamma'\gamma}(p)} \,\,\circ\\
		 &\quad \circ \,\rho(\gamma')(\tra_\gamma(p)):1_{F\circ \tra_{\gamma'\gamma}(p)}
		\end{align}
		The first term on the right hand side can be rewritten to
		\mathindent5mm
		\begin{equation*}
			1_{\tra_{\gamma'}F\tra_\gamma(p):F\tra_{\gamma'\gamma}(p)}\cdot \left(\rho(\gamma)(p):1_{F\circ\tra_{\gamma}(p)}\right) = 1_{t\left(\rho(\gamma')(\tra_\gamma(p)):1_{F\circ \tra_{\gamma'\gamma}(p)}\right)}\cdot \left(\rho(\gamma)(p):1_{F\circ\tra_{\gamma}(p)}\right).
		\end{equation*}\mathindent12.5mm
		Using the definition of composition in 2-groups in terms of crossed modules we note that for $h,h'\in H$:
		\begin{equation}\label{eq:comptomult}
			1_{t(h)}h'\circ h = h'\cdot h.
		\end{equation}
		Using this fact we finally obtain 
		\begin{equation}\label{eq:rhoHfunct}
			\rho_H(\gamma'\gamma)(p) = \rho_H(\gamma)(p)\cdot \rho_H(\gamma')(\tra_\gamma(p)).
		\end{equation}
		Compare this relation to Lemma \ref{lem:etaHcomposition}. There are several things to note about this relation. First of all the order of $\gamma$ and $\gamma'$ is inverted between left and right hand side, therefore to get something functorial we should use $\rho_H(\gamma)^{-1}$ instead. Since $\rho(\gamma)$ takes values on the principal bundle, it's more natural to work with paths on $\mc P$. Let $\pi\colon \mc P\to M$ be the projection, and let $\gamma\colon p\to q$, then we should consider $\rho(\pi_*\gamma)(p)^{-1}$. Suppose $\gamma\colon p\to q$ and $\gamma'\colon q\to r$, then note that in the second factor of the right hand side we have $\tra_\gamma(p)$ as argument, whereas we expect $q$ as argument for $\rho_H(\gamma')$. To remedy this, we instead consider something $G\ltimes H$-valued, with $\tra_{\pi_*\gamma}(p):q$ in the $G$ factor. We define for $\gamma\colon p\to q$:
		\begin{equation}
			\rho_{G\ltimes H}(\gamma) = \left(q:\tra_{\pi_*\gamma}(p),\,\rho_H(\pi_*\gamma)(p)\right)^{-1}=\left(\tra_{\pi_*\gamma}(p):q,\,\rho_H(\pi_*\gamma)\left(\tra_{\pi_*\gamma}^{-1}(q)\right)\right).
		\end{equation}
		Where we used the fact that $\alpha_{g^{-1}}\rho_H(\gamma)(p) = \rho_H(\gamma)(p\cdot g)$ by lemma \ref{lem:g2tor2morphH}. Now we check that this is functorial:
		\begin{align*}
			\rho_{G\ltimes H}(\gamma')\rho_{G\ltimes H}(\gamma) &= \left(r:\tra_{\pi_*\gamma'}(q),\,\rho_H(\pi_*\gamma')(q)\right)^{-1}\left(q:\tra_{\pi_*\gamma}(p),\,\rho_H(\pi_*\gamma)(p)\right)^{-1}\\
			&=\left((r:\tra_{\pi_*\gamma'}(q))\cdot (q:\tra_{\pi_*\gamma}(p)),\, \rho_H(\pi_*\gamma)(p)\cdot \alpha_{q:\tra_{\pi_*\gamma}(p)}\rho_H(\pi_*\gamma')(q)\right)^{-1}\\
			&=\left(r:\tra_{\pi_*\gamma'\gamma(q)},\,   \rho_H(\pi_*\gamma)(p)\rho_H(\pi_*\gamma')(\tra_{\pi_*\gamma}(p))\right)^{-1}\\
			&=\rho_{G\ltimes H}(\gamma'\gamma).
		\end{align*}
		Therefore $\rho_{G\ltimes H}$ defines a (smooth) functor $\Pi_1^\thin(\mc P_0)\to B (G{\ltimes} H)$. Thus $\rho_{G\ltimes H}$ can be equivalently described by a 1-form $\psi\in \Omega^1(\mc P_0,\mf g\ltimes \mf h)$. Such a form is determined by its $\mf g$-valued component and its $\mf h$-valued component. Recall that
		\[
			A(\gamma'(0)) = \left.\frac{\dd}{\dd t}\right|_{t=0} \tra_{\pi_*\gamma_t}^{-1}(\gamma(t)):\gamma(0)=\left.\frac{\dd}{\dd t}\right|_{t=0} \gamma(t):\tra_{\pi_*\gamma_t}(\gamma(0)).
		\]
		Therefore the $\mf g$ component of $\psi$ is just the connection $A$. The $\mf h$ component is the form $\phi\in \Omega^1(\mc P_0,\mf h)$ needed for this theorem:
		\begin{equation}\label{eq:defPhi}
			\phi(\gamma'(0)) := -\left.\frac{\dd}{\dd t}\right|_{t=0}\rho_H(\pi_*\gamma_t)(\gamma(0)).
		\end{equation}
		The sign of $\phi$ is purely convention, and is chosen to give a nicer formula for the gauge transform of $B$. 	The point of the computation so far is that we can also go back: given a form $\phi\in \Omega^1(\mc P_0, \mf h)$ we can construct a 1-morphism of transport functors, where the natural transformation $\rho(\gamma)$ is obtained for each $\gamma$ by taking a path-ordered exponential of $\phi$. Moreover this $\rho(\gamma)$ is compatible with composition of paths in the sense of relation \eqref{eq:1morphcomposition}. Hence $\rho(\gamma)$ defines a 1-morphism of transport 2-functors if we also confirm equation \ref{eq:1morph2morph}.
	
		The next step is to prove in which way $(F,\phi)$ transforms the connection $(A,B)$. This is achieved by differentiating relation \eqref{eq:defrhogamma} defining $t(\rho_H(\gamma)(p))$:
		\begin{equation}
			 F(\tra_A(\gamma)(p))\cdot t(\rho_H(\gamma)(p)) = \tra_{A'}(\gamma)(F(p)).
		\end{equation}
		Equivalently dividing by $F(q)=F(\gamma(1))$ we get
		\begin{equation}
			 \left(\tra_A(\gamma)(p):q\right) \cdot t(\rho_H(\gamma)(p)) = \tra_{A'}(\gamma)(F(p)):F(q).
		\end{equation}
		Then replacing $\gamma$ by $\gamma_t$ and differentiating we get
		\begin{equation}
			-A-t_*\phi=-F^*A'.
		\end{equation}
		Conversely any $\phi$ satisfying this relation will define a natural transformation $\rho(\gamma):F\circ \tra_A(\gamma)\to \tra_{A'}(\gamma)\circ F$. Finally we will show that the claimed relation between $B$ and $B'$ is equivalent to the compatibility of $\rho(\gamma)$ with 2-morphisms as in \eqref{eq:1morph2morph}. That is, for any bigon $\Sigma\colon \gamma\Rightarrow\gamma'$ we have the following relation:
	
		\begin{equation}
			{\begin{tikzcd}[column sep=8em, row sep =4em]
				\mc P_x \ar["\tra^A_\gamma"]{r} \ar["F"']{d} & \mc P_y \ar["F"]{d} \ar["\rho(\gamma)" description,shorten >=0.3cm,shorten <=0.3cm, Rightarrow]{dl}\\
				\mc P'_x \ar["\tra^{A'}_\gamma"{name=U,description}]{r} \ar[bend right=80, "\tra^{A'}_{\gamma'}"{name=D,below}]{r} \ar[Rightarrow, shorten >=0.1cm, from=U, to=D,"\tra^{2}_{(A',B')}\Sigma"{description, pos=0.4}]& \mc P_y'
			\end{tikzcd}}
			\quad=\quad
			{\begin{tikzcd}[column sep=8em, row sep =4em]
				\mc P_x \ar["\tra^A_{\gamma'}"{name=D,description}]{r} \ar["F"']{d} \ar[bend left = 80,"\tra^A_\gamma"{name=U}]{r} \ar[Rightarrow, from=U, to=D,shorten <=0.1cm, "\tra^2_{(A,B)}\Sigma"{description,pos=0.5}]& \mc P_y \ar["F"]{d} \ar["\rho(\gamma')" description,shorten >=0.3cm,shorten <=0.3cm, Rightarrow]{dl}\\
				\mc P_x' \ar["\tra^{A'}_{\gamma'}"']{r} & \mc P_y'.
			\end{tikzcd}}
		\end{equation}
		
		\noindent In equations this equality is
		\begin{equation}\label{eq:B2morphcompat}
			\begin{aligned}
				\tra^2_{(A',B')}\Sigma(F(p))\circ \rho(\gamma)(p) &= \rho(\gamma')(p)\circ F\left(\tra^2_{(A,B)}\Sigma(p)\right)\\
				&=\rho(\gamma')(p)\circ 1_{F\circ \tra_A\gamma(p)}\cdot \tra^2_{(A,B)}\Sigma(p):1_{ \tra_A\gamma(p)}
			\end{aligned}
		\end{equation}
		Then dividing both sides by the source $1_{F\tra_\gamma^A(p)}$ we obtain for the left side:
		\begin{align*}
			&\left(\tra^2_{(A',B')}\Sigma(F(p))\circ\rho(\gamma)(p)\right):1_{F\tra_\gamma^A(p)}\\
		 &\qquad = \left(\tra^2_{(A',B')}\Sigma(F(p)):1_{F\tra_\gamma^A(p)}\right)\circ \rho_H(\gamma)(p)\\
		 &\qquad = \left(1_{\tra^{A'}_\gamma(F(p)):F\tra_\gamma^A(p)}\cdot \left(\tra^2_{(A',B')}\Sigma(F(p)):1_{\tra^{A'}_\gamma(F(p))} \right)\right) \circ \rho_H(\gamma)(p)\\
		 &\qquad =  \left(\tra^2_{(A',B')}\Sigma(F(p)):1_{\tra^{A'}_\gamma(F(p))} \right) \cdot  \rho_H(\gamma)(p),
		\end{align*}
		where in the last step we used the relation \eqref{eq:comptomult}. Note that the first term on the last line lies in $H$. Using the same trick we obtain a similar expression for the right hand side of \eqref{eq:B2morphcompat}. Putting this together gives:
		\begin{equation}\label{eq:B2morphcompat2}
			\left(\tra^2_{(A',B')}\Sigma(F(p)):1_{\tra^{A'}_\gamma(F(p))} \right) \cdot  \rho_H(\gamma)(p) = \rho_H(\gamma')(p)\cdot \left(\tra^2_{(A,B)}\Sigma(p):1_{\tra^A_\gamma(p)}\right).
		\end{equation}
		To shorten notation let us denote
		\begin{equation}
			h_\Sigma:=\left(\tra^2_{(A,B)}\pi_*\Sigma(p):1_{\tra^{A}_{\pi_*\gamma}(p)} \right),\qquad h'_\Sigma = \left(\tra^2_{(A',B')}\pi_*\Sigma(F(p)):1_{\tra^{A'}_{\pi_*\gamma}(F(p))} \right).
		\end{equation}
		Furthermore let us also denote:
		\begin{equation}\label{eq:ggammadef}
			g_\gamma:= \tra_{\pi_*\gamma}^A(\gamma(0)):\gamma(1).
		\end{equation}
		Then finally let us denote
		\begin{equation}
			\hat \rho_H(\gamma) = \rho_H(\pi_*\gamma)(\gamma(0)).
		\end{equation}
		In this notation the functoriality of $\rho_H$ (cf. eq. \eqref{eq:rhoHfunct}) can be written as:
		\begin{equation}\label{eq:rhohatHfunct}
			\hat \rho_H(\gamma'\gamma)=\hat \rho_H(\gamma)\alpha_{g_\gamma^{-1}}\hat\rho_H(\gamma').
		\end{equation}
		
		\begin{figure}[!htb]
			\centering
			\includegraphics[width = 0.3\textwidth]{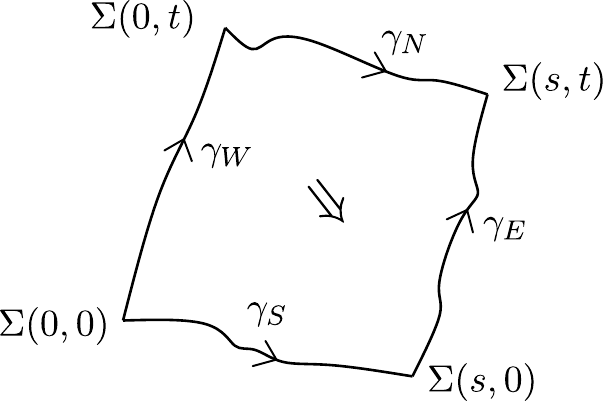}
			\caption{\label{fig:gammaNESW} The source an target of the bigon $\Sigma_*\Gamma(s,t)$ can be naturally split in two, and we label the four resulting paths by the cardinal directions.}
		\end{figure}
		
		Let $\Sigma$ be a smooth map $I^2\to P$, and consider the bigon $\Sigma_*\Gamma(s,t)$, with $\Gamma(s,t)$ as in figure \ref{fig:canonicalBigon}. Then denote the four boundary paths of this bigon by $\gamma_E,\gamma_N,\gamma_W,\gamma_S$ (for East, North, West and South) as sketched in figure \ref{fig:gammaNESW}.	Then equation \eqref{eq:B2morphcompat2} gives
		\begin{equation}
			\hat\rho_H(\gamma_E\gamma_S)\cdot h_{\Sigma_*\Gamma(s,t)}=h'_{\Sigma_*\Gamma(s,t)}\cdot \hat\rho_H(\gamma_N\gamma_W).
		\end{equation}
		Using the functoriality \eqref{eq:rhohatHfunct} for $\rho_H$ we can write this as 
		\begin{equation}\label{eq:rhosigmast}
			\hat \rho_H(\gamma_S)\alpha\left(g(\gamma_S)^{-1},\hat\rho_H(\gamma_E)\right)h_{\Sigma_*\Gamma(s,t)}=
			h'_{\Sigma_*\Gamma(s,t)}\hat\rho_H(\gamma_W)\alpha\left(g(\gamma_W)^{-1},\hat \rho_H(\gamma_N)\right).
		\end{equation}
		Let us refer to the left hand side and right hand side of \eqref{eq:rhosigmast} as $\textbf{LHS}$ and $\textbf{RHS}$ respectively. We wish to compute the second derivative of both sides at the origin. Let $X=\Sigma_*\partial s$ and $Y=\Sigma_*\partial t$, then by equation \eqref{eq:defPhi} we have:
		\begin{align*}
			\left.\frac{\partial}{\partial s}\right|_0&\hat \rho_H(\gamma_S) = -\phi_{\Sigma(0,0)}(X),
			&&\hspace{-0.24ex}\left.\frac{\partial}{\partial s}\right|_0\hat \rho_H(\gamma_N) = -\phi_{\Sigma(0,t)}(X),\\
			\left.\frac{\partial}{\partial t}\right|_0&\hat \rho_H(\gamma_W) = -\phi_{\Sigma(0,0)}(Y),
			&&\left.\frac{\partial}{\partial t}\right|_0\hat \rho_H(\gamma_E) = -\phi_{\Sigma(s,0)}(Y),\\
			\left.\frac{\partial}{\partial s}\right|_0& g(\gamma_S)^{-1} = A_{\Sigma(0,0)}(X),
			&&\left.\frac{\partial}{\partial t}\right|_0 g(\gamma_W)^{-1} = A_{\Sigma(0,0)}(Y).
		\end{align*}
		\noindent Let us begin with the left hand side, and compute first the derivative with respect to $s$.
		\begin{align*}
			\frac{\partial \textbf{LHS}}{\partial s} &= \frac{\partial\hat \rho_H(\gamma_S)}{\partial s}\alpha\left(g(\gamma_S)^{-1},\hat\rho_H(\gamma_E)\right)h_{\Sigma_*\Gamma(s,t)}\\
			&+
			\hat \rho_H(\gamma_S)\left(\alpha\left(\frac{\partial g(\gamma_S)^{-1}}{\partial s},\hat\rho_H(\gamma_E)\right) +\alpha\left(g(\gamma_S)^{-1},\frac{\partial \hat\rho_H(\gamma_E)}{\partial s}\right)\right) h_{\Sigma_*\Gamma(s,t)}\\
			&+\hat\rho_H(\gamma_S)\alpha\left(g(\gamma_S)^{-1},\hat\rho_H(\gamma_E)\right)\frac{\partial h_{\Sigma_*\Gamma(s,t)}}{\partial s}. 
		\end{align*}
		Consequently we can take the second derivative at 0.
		\begin{equation}
			\left. \frac{\partial^2 \textbf{LHS}}{\partial s\partial t} \right|_{(0,0)} = -\phi(X)\phi(Y)-\alpha_*(A(X),\phi(Y))-\frac{\partial \phi(Y)}{\partial s}-B(X,Y).
		\end{equation}
		Here we used the definition of $A$, $B$ and $\phi$ in terms of parallel transport. We will now compute the first derivative of the right hand side:
		\begin{align*}
			\frac{\partial \textbf{RHS}}{\partial s} &=
			\frac{\partial h'_{\Sigma_*\Gamma(s,t)}}{\partial s}\hat\rho_H(\gamma_W)\alpha\left(g(\gamma_W)^{-1},\hat \rho_H(\gamma_N)\right)\\
			&+h'_{\Sigma_*\Gamma(s,t)}\frac{\partial \hat\rho_H(\gamma_W)}{\partial s} \alpha\left(g(\gamma_W)^{-1},\hat \rho_H(\gamma_N)\right)\\
			&+h'_{\Sigma_*\Gamma(s,t)}\hat\rho_H(\gamma_W) \left(
			\alpha\left(\frac{\partial g(\gamma_W)^{-1}}{\partial s},\hat \rho_H(\gamma_N)\right)+
			\alpha\left(g(\gamma_W)^{-1},\frac{ \partial \hat \rho_H(\gamma_N)}{\partial s}\right)
			\right).
		\end{align*}
		Now let us compute the second derivative:
		\begin{equation}
			\left.\frac{\partial^2 \textbf{RHS}}{\partial s\partial t}\right|_{(0,0)} = -F^*B'(X,Y)-\frac{\partial \phi(X)}{\partial t}-\phi(Y)\phi(X)-\alpha_*(A(Y),\phi(X)).
		\end{equation}
		Now combining the two this gives:
		\begin{equation}
			\begin{aligned}
				F^*B'(X,Y) &= B(X,Y)+\left(\frac{\partial \phi(Y)}{\partial s}-\frac{\partial \phi(X)}{\partial t}\right) +\left(\phi(X)\phi(Y)-\phi(Y)\phi(X)\right)+\\ &+\left(\alpha_*(A(X),\phi(Y))-\alpha_*(A(Y),\phi(X))\right).\\
			\end{aligned}
		\end{equation}
		Recall equation \eqref{eq:infcrossedmodulebracket} defining the Lie bracket of $\mf g\ltimes \mf h$, using this we identify the last term with ${[A\wedge \phi](X,Y)}$. Since this holds for all smooth $\Gamma:I^2\to P$ we conclude
		\begin{equation*}
			F^*B' = B+d\phi+\frac 12[\phi,\phi]+\frac12[A\wedge \phi].
		\end{equation*}

		Thus we have shown that given a 1-morphism of transport 2-functors we naturally obtain a 1-morphism of the corresponding principal $\mc G$-2-bundle with 2-connection, i.e. we constructed the action of a 2-functor $\Trans^2(M,\mc G)\to\Bun_\nabla^2(M,\mc G)$ on the 1-morphisms. We also showed that through an integration procedure we obtain a 2-functor in the opposite direction, except for showing that equation \eqref{eq:1morph2morph} holds. We claim that the arguments provided in \cite{SW11} can be adapted with only minor changes to show this. \qed

		Let $\ms F$ and $\ms G$ be transport 2-functors $\Pi_2^\thin(M)\to \mc G\tor$ corresponding to principal $\mc G$-2-bundles with connection $(\mc P,A, B)$ and $(\mc P', A', B')$ respectively. Consider a 2-morphism
		\begin{equation}
			\begin{tikzcd}[column sep=5em,row sep=6.5em]
				\ms F \ar["{\rho}"{name=U},bend left= 50]{r} \ar[bend right = 50, "{\rho'}"{name = D,below}]{r} \ar[from= U,to=D,Rightarrow,"\mc A"{description},shorten <=1mm,shorten >=1mm] & \ms G
			\end{tikzcd}\quad \leftrightarrow \quad 
			\begin{tikzcd}[column sep=2em]
				(\mc P, A, B) \ar["{(F,\phi)}"{name=U},bend left= 50]{r} \ar[bend right = 50, "{(F',\phi')}"{name = D,below}]{r} \ar[from= U,to=D,Rightarrow,"a"{description},shorten <=1mm,shorten >=1mm] & (\mc P',A',B')
			\end{tikzcd}
		\end{equation}
		Such a 2-morphism $\mc A\colon \rho\rightarrow\rho'$ (cf. lem. \ref{lem:g2tor2morphH}) can be equivalently described by a $G$-equivariant map $a\colon \mc P_0\to H$ (i.e. $a(p\cdot g) = \alpha_{g^{-1}}a(p)$).  such that
		\begin{equation}
			\phi'=\Ad_a \phi-R_{a^{-1}}\circ(\alpha_a)_*A+\dd a a^{-1},
		\end{equation}
		where $(\alpha_a)_*\colon \mf g\to \mf h$ is obtained from $\alpha\colon G\times H\to H$ by fixing the second coordinate and differentiating.
	
		\begin{lemma}\label{lem:2functor2morphism}
			If $\mc A\colon \rho\rightarrow \rho'$ is a 2-morphism of transport 2-functors, such that $\rho$ and $\rho'$ correspond to $(F,\phi)$\, $(F',\phi')$ by lemma \ref{lem:2functor1morphism} then $a\colon \mc P_0\to H$ satisfies
			\begin{equation}\label{eq:ais2morphism}
				\phi'=\Ad_a \phi-R_{a^{-1}}\circ(\alpha_a)_*A+\dd a a^{-1},
			\end{equation} 
			and hence defines a 2-morphism of principal 2-bundles with 2-connection. Conversely if $a$ satisfies \eqref{eq:ais2morphism} then $\mc A$ is a 2-morphism of transport 2-functors. 
		\end{lemma}
	
		\proof Recall from definition \ref{def:2trans2morph} that $\mc A$ assigns to each $x\in M$ a natural transformation $\ms F_x\rightarrow \ms G_x$, that is for each $p\in \mc P_{0,x}$:
		\[
			\mc A(x)(p)\colon \ms F(p)\to \ms G(p).
		\]
		By Lemma \ref{lem:g2tor2morphH} $\mc A$ is equivalent to the equivariant map:
		\begin{equation}
			a\colon \mc P_0\to H,\qquad a(p)=\mc A(\pi(p))(p):1_{\ms F(p)}.
		\end{equation}
		For any $\gamma\colon x\to y$ we then have the following relation (cf. eq. \eqref{eq:2morphtransformrho}):
		\begin{equation}
			\begin{tikzcd}[column sep=6em, row sep =6em]
				\mc P_x \ar["\tra^A_\gamma"]{r} \ar["\ms G"{left,name=L},bend right=100]{d} \ar["\ms F"{description,name=R}]{d} \ar[from =R, to=L, Rightarrow, "\mc A(x)"{above,pos=0.4}, shorten >= 0.1cm] & \mc P_y \ar["\ms F"]{d} \ar["\rho(\gamma)" description,shorten >=0.3cm,shorten <=0.3cm, Rightarrow]{dl}\\
				\mc P_x' \ar["\tra_\gamma^{A'}"']{r} & \mc P_y'
			\end{tikzcd}
			\quad=\quad
			\begin{tikzcd}[column sep=6em, row sep =6em]
				\mc P_x \ar["\tra_\gamma^A"]{r} \ar["\ms G"']{d} & \mc P_y \ar["\ms G"{description, name = L}]{d} \ar["\ms F"{right, name=R}, bend left = 100]{d} \ar["\mc A(y)"{above,pos=0.5},from=R, to=L,Rightarrow,shorten <=0.1cm] \ar["\rho'(\gamma)" description,shorten >=0.3cm,shorten <=0.3cm, Rightarrow]{dl}\\
				\mc P_x' \ar["\tra_\gamma^{A'}"']{r} & \mc P_y'
			\end{tikzcd}
		\end{equation}
		Writing this out we obtain for $p\in \mc P_{x,0}$:
		\begin{equation}
			\left(\tra^{A'}_\gamma\!\circ\, \mc A(x)\right)\circ \rho(\gamma)(p) = \rho'(\gamma)\circ \mc A(y)(\tra_\gamma^A(p)).
		\end{equation}
		To obtain a more useful equation we will divide both sides by $1_{\ms F\circ\tra^A_\gamma(p)}$ and simplify the resulting expression. For the left hand side we obtain:
		\begin{align*}
			1_{\tra^{A'}_\gamma(\ms F(p)):\ms F\circ \tra^A_\gamma(p)}\cdot  a(p) \circ \rho_H(\gamma)(p) = a(p)\cdot \rho_H(\gamma)(p).
		\end{align*}
		For the right hand side we obtain:
		\begin{align*}
			\rho'(\gamma)(p):1_{\ms F\circ \tra^A_\gamma(p)}\circ a(\tra^A_\gamma(p)) &= \rho_H'(\gamma)(p)\cdot 1_{\ms G\circ \tra^A_\gamma(p):\ms F\circ \tra^A_\gamma(p)}\circ a(\tra^A_\gamma(p))\\
			&= \rho_H'(\gamma)(p)\cdot a(\tra^A_\gamma(p)).
		\end{align*}
		If instead we consider a path $\gamma_t:p\to \gamma(t)$, then this gives equality:
		\begin{equation}
			\alpha(g(\gamma_t)^{-1},\,\rho_H'(\pi_*\gamma_t)(p))=\alpha\Big(g(\gamma_t)^{-1},\,a(p)\cdot\rho_H(\pi_*\gamma_t)(p)\Big)a(\gamma(t))^{-1},
		\end{equation}
		where $g_\gamma$ is defined in the proof of the previous theorem (eq. \eqref{eq:ggammadef}). The equation is arranged in this way such that both sides are $1$ at $t=0$. Differentiating this equation at 0 gives:
		\begin{equation}
			\phi'=\Ad_a\phi-r_{a^{-1}}(\alpha_a)_*A+daa^{-1},
		\end{equation}
		w hich is exactly what we wanted to prove. We naturally constructed a 2-morphism of principal $\mc G$-2-bundles with 2-connection out of a 2-morphism of transport 2-functors. That is, we have a 2-functor $\Trans^2(M,\mc G)\to\Bun_\nabla^2(M,\mc G)$. To construct a 2-functor in the opposite direction we again claim that the arguments provided in \cite{SW11} can be adapted with only minor changes. \qed

		Lemmas \ref{lem:2functor1morphism} and \ref{lem:2functor2morphism} now show that the equivalence $\Bun_\nabla^1(M,G)\cong \Trans^1(M,G)$ of theorem \ref{thm:transportisbundle} extends to 2-functors $\Trans^2(M,\mc G)\to\Bun_\nabla^2(M,\mc G)$ and $\Bun_\nabla^2(M,\mc G)\to \Trans^2(M,\mc G)$. Since these functors are inverse to each other on the level of 1- and 2-morphisms, this gives an equivalence of 2-categories.
		
		\begin{theorem}
			There is an equivalence of 2-categories $\Trans^2(M,\mc G)\to\Bun_\nabla^2(M,\mc G)$.
		\end{theorem}
		
		This can likely be upgraded to an equivalence of (2-)stacks $\Trans^2(\mc G)\to \Bun_\nabla^2(\mc G)$, however we do not know how to prove this. Furthermore one can try replacing the strict 2-groupoid $\Pi_2^\thin(M)$ with a weak version (by not identifying thinly homotopic paths), and with some modifications this should give a weaker notion of parallel 2-transport.  One can also make the notion of $\mc G$-torsor weaker by posing that the functor $(\pi_1\times  R)\colon \mc X\times G\to \mc X\times \mc X$ is an equivalence (and not invertible on the nose). This is investigated in \cite{Wal16,Wal17}. Both these modifications are more difficult to understand because then there is no transport functor underlying a 2-transport functor, or equivalently a principal 2-bundle does not have an associated principal bundle. One can also change the target category entirely, and use a notion of local trivializations of 2-functors as in \cite{SW13}. It is likely that our notion of 2-transport functor coincides with that of \cite{SW13}, if we specify the target category to be $\mc G\tor$. Finally one can also try to adapt the definition of a 2-bundle to suit weak 2-groups.
		
	\section{Curvature and higher non-Abelian Stokes}
		We report a result obtained by the author in \cite{Voorhaar}. Given the importance of the non-Abelian Stokes Theorem one can ask if there is a higher version of this. This is not a new result, and can be found for instance in \cite{MartinsPicken1, MartinsPicken2, AFG} in different context. Our result is essentially a corollary to lemma 2.16 in \cite{SW11}, suitably adapted to this setting. In this framework it allows us to prove a higher version of the Ambrose Singer theorem, and the non-Abelian Stokes Theorem is essential in generalizing 2-transport to $n$-transport.
		
		\begin{definition}
			The curvature $\tilde K\in \Omega^3(P,\mf h)$ of a 2-connection $(A,B)$ on a principal $\mc G$-2-bundle $\mc P$ is given by
			\begin{equation}
				\tilde K = dB+\frac 12[A\wedge B].
			\end{equation}
			Alternatively, denoting $h:T\mc P_0\to T\mc P_0$ the projection onto $\ker A\subset T\mc P_0$, it is given by $\tilde K=dB\circ h$.
		\end{definition}
		
		Note that $t_*\tilde K=dF+\frac12[A\wedge F]=0$ by the Bianchi identity, therefore $\tilde K$ takes values in $\ker t_*$.  Since $\ker t_*$ is an Abelian Lie algebra, $\tilde K$ is in fact basic. That is, there is a form $K\in \Omega^3(M,\ker t_*)$ such that $\tilde K=\pi^* K$. The higher non-Abelian Stokes Theorem then simply states that the difference between the 2-transport over two homotpic bigons is given by an integral over $K$. 
	
		\begin{theorem}[\textbf{higher non-Abelian Stokes}]\label{thm:nonabelianstokes2}
			Let $h:I^3\to M$ be a homotopy $h_0\Rrightarrow h_1$ of bigons and let $p\in P_x$ then
			\begin{equation}
				\tra^2(h_1)(p):\tra^2(h_0)(p)=\exp\int_{I^3}-h^* K
			\end{equation}
		\end{theorem}
		
		We remark that $\tra^2(h_i)(p)$ both share the same source and target, therefore the left hand side indeed lands in $\ker t$. Furthermore the right hand side involves an ordinary integral instead of a path ordered exponential, this is because $ K$ is Abelian. The proof proceeds similarly to our proof of the `classical' non-Abelian Stokes theorem \ref{thm:nonabelianstokes}. Most of the work in this proof involves translating the proof of lemma 2.16 in \cite{SW11} to our setting.
		
		\proof Let $h_u$ denote the bigon $h(u,\cdot,\cdot)$. Then we define:
		\[
			f(\rho) = \tra^2(h_0)(p):\tra^2(h_\rho)(p).
		\]
		This function satisfies
		\[
			f(\rho+r)f(\rho)^{-1} = \tra^2(h_{\rho})(p):\tra^2(h_{\rho+r})(p).
		\]
		And hence we obtain the differential equation
		\[
			f'(\rho) = \left[\left.\frac{\partial}{\partial r}\right|_0\tra^2(h_{\rho})(p):\tra^2(h_{\rho+r})(p)\right]\cdot f(\rho).
		\]
		Since the factors involved are abelian, we have
		\[
			f(1)^{-1} = \exp\int_0^1\!\dd \rho\,-\left.\frac{\partial}{\partial r}\right|_0\tra^2(h_{\rho})(p):\tra^2(h_{\rho+r})(p).
		\]
		\begin{figure}[htb]
		\centering
		\includegraphics[width=0.7\textwidth]{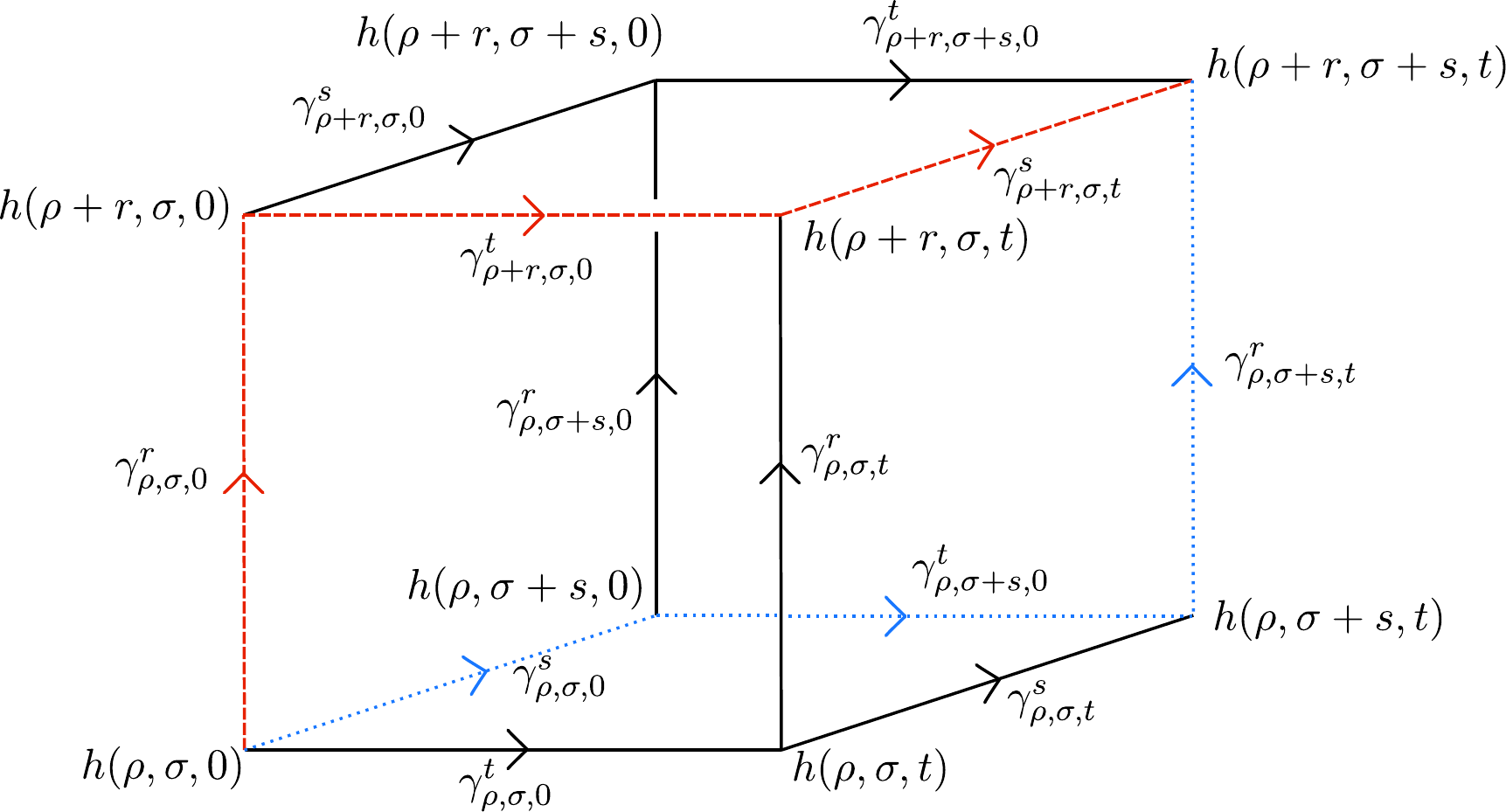}
		\caption{\label{fig:pathCube} Diagram of the image of a cube in $I^3$ under a map $h:I^3\to M$. The paths forming the edges of the cube are labeled. The boundary of the cube is split into two bigons, the source and target of which are indicated by blue and red on the diagram.}
		\end{figure}
		\begin{figure}[htb]
		\centering
		\includegraphics[width=0.8\textwidth]{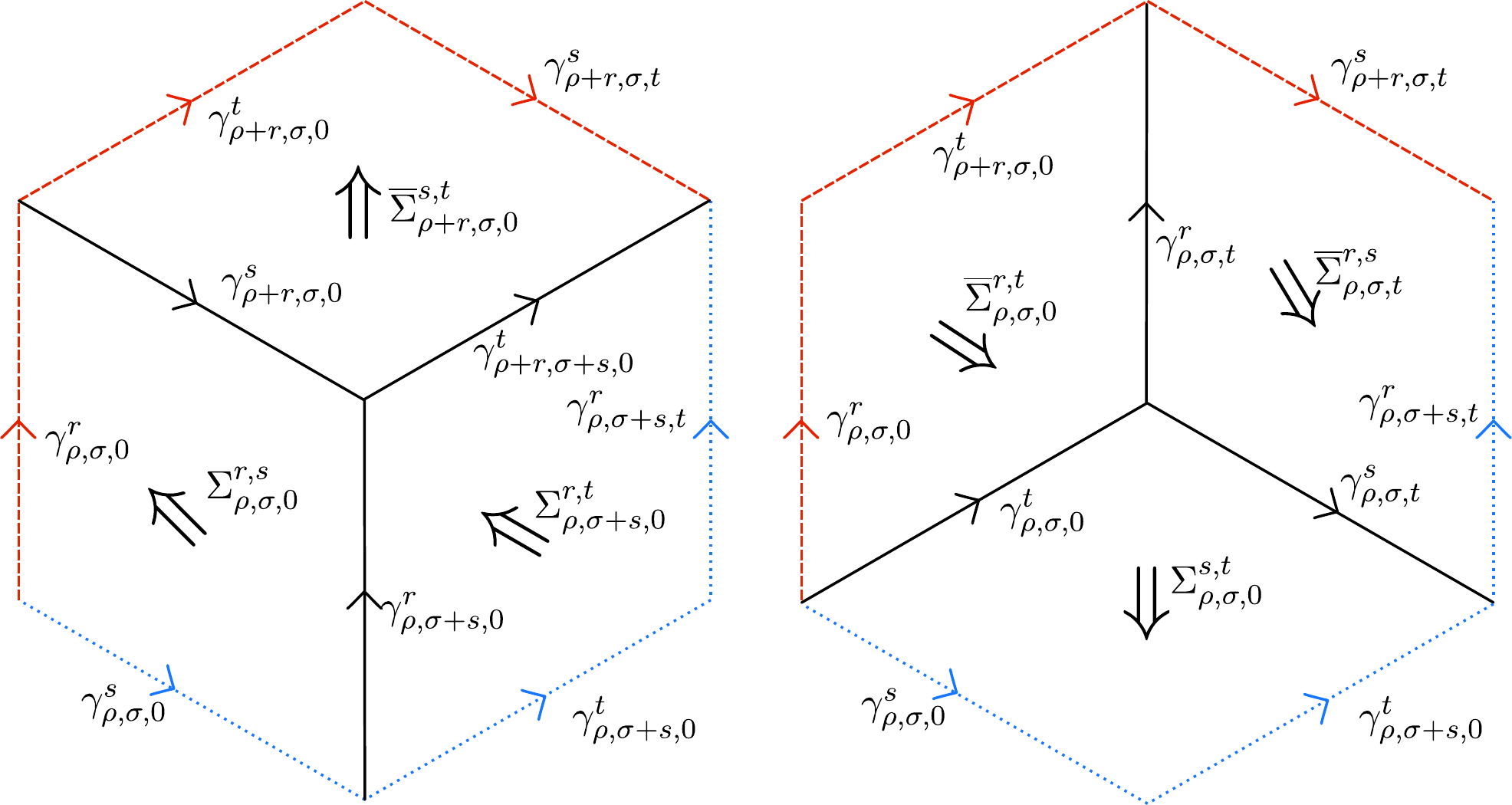}
		\caption{\label{fig:bigonCube} The two bigons forming the boundary of the cube in figure \ref{fig:pathCube}, each the composition of three bigons given by the faces of the cube.}
		\end{figure}
		
		The idea is now to add two more parameters, and take the holonomy around a small cube. Consider the cube drawn in figure \ref{fig:pathCube}. We associate a bigon to all the faces of the cube as shown in figure \ref{fig:bigonCube}. The two halves of the cube shown in figure \ref{fig:bigonCube} both represent a single bigon between the same paths. Let us call these bigons respectively $L(\rho,r,\sigma,s,t)$, $R(\rho,r,\sigma,s,t)$. Note that $L(\rho,r,0,1,1)=\overline{h_{\rho+r}}$ and $R(\rho,r,0,1,1) =  h_{\rho}$. 
		
		For any bigon $\Sigma:\gamma\Rightarrow\gamma'$ we denote 
		\begin{equation}
			\tra^2_H(\Sigma,p) = \tra^2(\Sigma)(p):\id_{\tra^1(\gamma)(p)}\in H.
		\end{equation}
		This notation comes from lemma \ref{lem:g2tor2morphH}, and we refer to lemma \ref{lem:etaHcomposition} for its vertical and horizontal composition rules. We are now interested in
		\begin{align}
			u_{\rho,\sigma,0}&(r,s,t):=\tra^2(L(\rho,r,\sigma,s,t))(p):\tra^2(\overline R(\rho,r,\sigma,s,t))(p)\\
			 &= \tra^2_H(L,p)\tra^2_H(R,p)\\
			 &= \tra^2_H(\Sigma{}^{r,t}_{\rho,\sigma+s,0},p) \tra^2_H(\Sigma^{s,t}_{\rho+r,\sigma,0},p)^{-1} \tra^2_H(\Sigma^{r,t}_{\rho,\sigma,0})^{-1}\\
			 &\qquad \cdot \tra^2_H(\Sigma^{r,s}_{\rho,\sigma,t},\tra^1(\gamma^t_{\rho,\sigma,0})(p))^{-1} \tra^2_H(\Sigma^{s,t}_{\rho,\sigma,0},p).
		\end{align}
		We will need the following two properties of this function:
		\begin{lemma}\label{lem:uproperties}
			The function $u_{\rho,\sigma,0}(r,s,t)$ satisfies:
		\begin{align} 
			\label{eq:uproperty1}\left.\frac{\partial}{\partial s}\right|_0u_{\rho,\sigma,0}(r,s,1) = \frac{\partial}{\partial \sigma}u_{\rho,0,0}(r,\sigma,1),\\
			\label{eq:uproperty2}\left.\frac{\partial^3}{\partial r\partial s \partial t}\right|_{0,0,0}u_{\rho,\sigma,0}(r,s,\tau+t) = \tilde h^*\tilde K_{\rho,\sigma,\tau},
		\end{align}
		where $\tilde h(\rho,\sigma,\tau) = \tra^1(\gamma^\tau_{\rho,\sigma,0})(p)$.
		\end{lemma}
		
		Assuming this lemma we prove the theorem:
		
		\begin{align*}
			f(1)^{-1} &= \exp\int_0^1\!\!\dd\rho\,-\left.\frac{\partial}{\partial r}\right|_0 u_{\rho,0,0}(r,1,1)\\
				 &= \exp\int_0^1\!\!\dd\rho\int_0^1\!\!\dd \sigma\,-\left.\frac{\partial}{\partial r}\right|_0 \frac{\partial}{\partial \sigma}u_{\rho,0,0}(r,\sigma,1)\\
				 &= \exp\int_0^1\!\!\dd\rho\int_0^1\!\!\dd \sigma\,-\left.\frac{\partial^2}{\partial r\partial s}\right|_{0,0}u_{\rho,\sigma,0}(r,s,1)\\
				 &= \exp\int_0^1\!\!\dd\rho\int_0^1\!\!\dd \sigma\int_0^1\!\!\dd\tau\,- \left.\frac{\partial^3}{\partial r\partial s\partial t}\right|_{0,0,0}u_{\rho,\sigma,0}(r,s,\tau+t)\\
				 &= \exp\int_0^1\!\!\dd\rho\int_0^1\!\!\dd \sigma\int_0^1\!\!\dd\tau\, \,-\tilde h^*\tilde K = \exp\int_{I^3}\!\!-h^* K.\qedhere
		\end{align*}
		
		\proof[Proof of lemma \ref{lem:uproperties}] We first prove equation \eqref{eq:uproperty1}. We compute:
		\begin{align*}
		u_{\rho,\sigma,0}(r,s,1)&=\tra^2_H(\Sigma{}^{r,1}_{\rho,\sigma+s,0},p) \tra^2_H(\Sigma^{s,1}_{\rho+r,\sigma,0},p)^{-1} \tra^2_H(\Sigma^{r,1}_{\rho,\sigma,0})^{-1}\\
					 &\qquad\cdot \tra^2_H(\Sigma^{r,s}_{\rho,\sigma,1},\tra^1(\gamma^1_{\rho,\sigma,0})(p))^{-1} \tra^2_H(\Sigma^{s,1}_{\rho,\sigma,0},p).
		\end{align*}
		Note that $\Sigma^{r,s}_{\rho,\sigma,1}$ is $\id_{h(-,-,1)}$, and hence does not contribute. Taking the derivative with respect to $s$ we obtain:
		\[
			\left.\frac{\partial}{\partial s}\right|_0u_{\rho,\sigma,0}(r,s,1) = \left.\frac{\partial}{\partial s}\right|_0\tra^2_H(\Sigma{}^{r,1}_{\rho,\sigma+s,0},p) \tra^2_H(\Sigma^{s,1}_{\rho+r,\sigma,0},p)^{-1} \tra^2_H(\Sigma^{s,1}_{\rho,\sigma,0},p).
		\]
		On the other hand we compute:
		\begin{align*}
			u_{\rho,0,0}(r,\sigma+s,1)&=\tra^2_H(\Sigma{}^{r,1}_{\rho,\sigma+s,0},p) \tra^2_H(\Sigma^{\sigma+s,1}_{\rho+r,0,0},p)^{-1} \tra^2_H(\Sigma^{r,1}_{\rho,0,0})^{-1}\\
						 &\qquad \cdot \tra^2_H(\Sigma^{r,\sigma+s}_{\rho,0,1},\tra^1(\gamma^1_{\rho,0,0})(p))^{-1} \tra^2_H(\Sigma^{\sigma+s,1}_{\rho,0,0},p).
		\end{align*}
		Again we note that $\Sigma^{r,\sigma +s}_{\rho,0,1}=\id_{h(-,-,1)}$, and furthermore we decompose $\Sigma^{\sigma+s,1}_{\rho,0,0}=\Sigma^{s,1}_{\rho,\sigma,0}\bullet\Sigma^{\sigma,1}_{\rho,0,0}$, and note that the right factor is independent of $s$. Now taking the derivative with respect to $s$, we obtain:
		\[
			\left.\frac{\partial}{\partial s}\right|_0u_{\rho,0,0}(r,\sigma+s,1) = \left.\frac{\partial}{\partial s}\right|_0\tra^2_H(\Sigma{}^{r,1}_{\rho,\sigma+s,0},p) \tra^2_H(\Sigma^{s,1}_{\rho+r,\sigma,0},p)^{-1} \tra^2_H(\Sigma^{s,1}_{\rho,\sigma,0},p).
		\]
		Which confirms equation \eqref{eq:uproperty1}. To prove equation \eqref{eq:uproperty2} we first compute $u_{\rho,\sigma,0}(r,s,\tau+t)$ and decompose all the bigons involving $t+\tau$ into two pieces:
		\begin{align}
			&\tra^2_H(\Sigma^{r,\tau+t}_{\rho,\sigma+s,0},p)&&\hspace{-8em}=\tra^2_H(\Sigma^{r,\tau}_{\rho,\sigma+s,0},p)\tra^2_H(\Sigma^{r,t}_{\rho,\sigma+s,\tau},\tra^1(\gamma^\tau_{\rho,\sigma+s,0})(p)),\\
			&\tra^2_H(\Sigma^{r,\tau+t}_{\rho,\sigma,0},p)^{-1}&&\hspace{-8em}=\tra^2_H(\Sigma^{r,t}_{\rho,\sigma,\tau},\tra^1(\gamma^\tau_{\rho,\sigma,0})(p))^{-1} \tra^2_H(\Sigma^{r,\tau}_{\rho,\sigma,0},p)^{-1},\\
			&\tra^2_H(\Sigma^{s,\tau+t}_{\rho,\sigma,0},p)&&\hspace{-8em}=\tra^2_H(\Sigma^{s,\tau}_{\rho,\sigma,0},p)\tra^2_H(\Sigma^{s,t}_{\rho,\sigma,\tau},\tra^1(\gamma^\tau_{\rho,\sigma,0})(p)),\\
			&\tra^2_H(\Sigma^{s,\tau+t}_{\rho+r,\sigma,0},p)^{-1}&&\hspace{-8em}=\tra^2_H(\Sigma^{s,t}_{\rho+r,\sigma,\tau},\tra^1(\gamma^\tau_{\rho+r,\sigma,0})(p))^{-1} \tra^2_H(\Sigma^{s,\tau}_{\rho+r,\sigma,0},p)^{-1}.
		\end{align}
		Thus we obtain:
		\begin{align*}
		u_{\rho,\sigma,0}(r,s,\tau+t) &= \tra^2_H(\Sigma^{r,\tau}_{\rho,\sigma+s,0},p)\tra^2_H(\Sigma^{r,t}_{\rho,\sigma+s,\tau},\tra^1(\gamma^\tau_{\rho,\sigma+s,0})(p))\\
		  	&\qquad \cdot \tra^2_H(\Sigma^{s,t}_{\rho+r,\sigma,\tau},\tra^1(\gamma^\tau_{\rho+r,\sigma,0})(p))^{-1} \tra^2_H(\Sigma^{s,\tau}_{\rho+r,\sigma,0},p)^{-1}\\
		  	&\qquad \cdot  \tra^2_H(\Sigma^{r,t}_{\rho,\sigma,\tau},\tra^1(\gamma^\tau_{\rho,\sigma,0})(p))^{-1} \tra^2_H(\Sigma^{r,\tau}_{\rho,\sigma,0},p)^{-1}\\
			&\qquad \cdot \tra^2_H(\Sigma^{r,s}_{\rho,\sigma,t+\tau},\tra^1(\gamma^{t+\tau}_{\rho,\sigma,0})(p))^{-1}\\
			&\qquad \cdot \tra^2_H(\Sigma^{s,\tau}_{\rho,\sigma,0},p)\tra^2_H(\Sigma^{s,t}_{\rho,\sigma,\tau},\tra^1(\gamma^\tau_{\rho,\sigma,0})(p)).
		\end{align*}
		To compute the third partial derivative of this expression we first note that for example,
		\begin{equation}
			1=\tra^2_H((\Sigma^{s,t}_{\rho,\sigma,0})^{-1}\circ \Sigma^{s,t}_{\rho,\sigma,0},p)=\tra^2_H((\Sigma^{s,t}_{\rho,\sigma,0})^{-1},\tilde h(\rho,\sigma+s,t))\tra^2_H(\Sigma^{s,t}_{\rho,\sigma,0},p).
		\end{equation}
		Therefore per definition of $B$ (cf. eq. \eqref{eq:Bdefinition}) we have
		\begin{equation}
			\left.\frac{\partial^2}{\partial s\partial t}\right|_{0,0}\tra^2_H(\Sigma^{s,t}_{\rho,\sigma,0},p) = -B_{\tilde h(\rho,\sigma,0)}(\tilde h_*\partial s,\tilde h^*\partial t).
		\end{equation}
		Using this we obtain
		\begin{align*}
			\left.\frac{\partial^3}{\partial r\partial s\partial t}\right|_{0,0,0}\tra^2_H(\Sigma^{r,t}_{\rho,\sigma+s,\tau},\tra^1(\gamma^\tau_{\rho,\sigma+s,0})(p))
			&=\left.\frac{\partial}{\partial s}\right|_0 -B_{\tilde h(\rho,\sigma+s,\tau)}(\tilde h_*\partial r,\tilde h_*\partial t)
			\\&=-\mc L_{\partial s}(\tilde h^*B)_{(\rho,\sigma,\tau)}(\partial r,\partial t).
		\end{align*}
		Identifying three similar terms we together obtain $(\tilde h^*dB)_{(\rho,\sigma,\tau)}(\partial r,\partial s,\partial t)$. Now there are two more terms contributing to the derivative. Consider
		\begin{align}
			&\left.\frac{\partial^3}{\partial r\partial s\partial t}\right|_{0,0,0}\tra^2_H(\Sigma^{r,\tau}_{\rho,\sigma+s,0},p) \tra^2_H(\Sigma^{s,t}_{\rho+r,\sigma,\tau},\tra^1(\gamma^\tau_{\rho+r,\sigma,0})(p))^{-1}\tra^2_H(\Sigma^{r,\tau}_{\rho,\sigma,0},p)^{-1}\\
			&=\left.\frac{\partial^3}{\partial r\partial s\partial t}\right|_{0,0,0}
				\alpha\left(t(\tra^2_H(\Sigma^{r,\tau}_{\rho,\sigma,0},p)),\,
				\tra^2_H(\Sigma^{s,t}_{\rho,\sigma,\tau},\tra^1(\gamma^\tau_{\rho,\sigma,0})(p))^{-1}\right)\\
			&=\left.\frac{\partial^3}{\partial r\partial s\partial t}\right|_{0,0,0}
				\alpha\left(\tilde h(\rho+r,\sigma,\tau):\tra(\gamma^r_{\rho,\sigma,\tau})(\tilde h(\rho,\sigma,\tau)),\,
				\tra^2_H(\Sigma^{s,t}_{\rho,\sigma,\tau},\tra^1(\gamma^\tau_{\rho,\sigma,0})(p))^{-1}\right)\\
			&=\alpha_*(A_{(\rho,\sigma,\tau)}(\tilde h^*\partial r),B_{(\rho,\sigma,\tau)}(\tilde h^*\partial s,\tilde h^*\partial t)).
		\end{align}
		Similarly we also obtain a factor
		\[
			-\alpha_*(A_{(\rho,\sigma,\tau)}(\tilde h^*\partial s),B_{(\rho,\sigma,\tau)}(\tilde h^*\partial r,\tilde h^*\partial t)),
		\]
		and we note that $\tilde h^*\partial t$ is a horizontal vector field, and is thus annihilated by $A$, meaning these two factors sum up to
		\[
			\alpha_*(\tilde h^*A,\tilde h^*B)(\partial r,\partial s,\partial t) = \frac12 [\tilde h^*A,\tilde h^*B](\partial r,\partial s,\partial t).
		\]
		Together with the $\tilde h^*dB$ factor, we obtain equation \eqref{eq:uproperty2}.\qed
		
		\begin{corollary}
			The 2-transport of a 2-connection $(A,B)$ as defined by equation \eqref{eq:2transportFormula} is thin-homotopy invariant. That is, if $\Sigma,\Sigma'$ are thinly homotopic then
			\[
				\tra^2(\Sigma)=\tra^2(\Sigma').
			\]
		\end{corollary}
	
		\proof Let $h:\Sigma\Rrightarrow \Sigma'$ be a rank 2 homotopy of bigons. Then $h^*K=0$, since $K$ is a rank 3 differential form. Hence
		\[
			\tra^2(\Sigma)(p)\tra^2(\Sigma')(p)=\exp\int_{I^3} -h^*K = 1.\qedhere
		\]
		
		\subsection{An Ambrose-Singer Theorem}
		
		Recall the classical Ambrose-Singer theorem \cite{KobayashiNomizu}.
		\begin{theorem}[Ambrose-Singer]
			Let $P\to M$ be a principal $G$-bundle with connection $A$, and let $p\in P_x$. Recall the definition of the \textit{holonomy group}:
			\begin{equation}
				\hol_p=\set[\tra_\gamma(p):p]{\gamma\in \Pi_1^\thin(M,x)}.
			\end{equation}
			Concatenation of paths defines a group operation. The Lie algebra $\mf{hol}_p$ of $\hol_p$ satisfies
			\begin{equation}
				\mf{hol}_p = \set[F_{\tra_\gamma(p)}(X,Y)]{\gamma\in \Pi_1^\thin(M),\gamma(0)=x,\,\,X,Y\in T_{\tra(\gamma,p)}P},
			\end{equation}
			where $F$ is the curvature of $A$. 
		\end{theorem}
		The identity component of $\hol_p$ is the \textit{reduced holonomy group} $\hat \hol_p$, i.e. the holonomies of all contractible groups. This theorem can then be rephrased as saying that the curvature completely determines the reduced holonomy group. This is not so surprising keeping the non-Abelian Stokes theorem \ref{thm:nonabelianstokes} in mind; the non-Abelian Stokes theorem tells us that the holonomy of a contractible loop can be expressed as an integral over the curvature. The precise argument showing how the Ambrose-Singer theorem is a corollary to the non-Abelian Stokes theorem appears in \cite{Voorhaar}. This provides a proof quite different from the original. The original prove relies on the fact that the curvature appears as the obstruction to involutivity of the horizontal distribution defined by a connection. There does not seem to be an analogous interpretation of the curvature of a 2-connection.
		
		Consequently we will shows that an Ambrose-Singer theorem for 2-gauge theory appears as a corollary to the higher non-Abelian Stokes Theorem. We first define the 2-holonomy group (n.b: this is a group, not a 2-group). 
		
		\begin{definition}
			Let $\mc P_1\rightrightarrows \mc P_0\to M$ be a principal $(G\ltimes H\rightrightarrows G)$-2-bundle with 2-connection $(A,B)$, and let $p\in (\mc P_0)_x$. Then we define the 2-holonomy group:
			\begin{equation}
				\hol^2_p=\set[\tra^2_H(\Sigma,p)]{\Sigma:\gamma\Rightarrow\gamma, \,\, s(\gamma)=x}\qquad \subset \ker t\subset Z(H),
			\end{equation}
			where we recall that $\tra^2_H(\Sigma,p)=\tra^2_\Sigma(p):\id_{\tra(s(\Sigma),p)}$. We also define the reduced 2-holonomy group $\hat{\hol^2_p}$ to be the subgroup of $\hol^2_p$ consisting of elements $\tra^2_H(\Sigma,p)$ where $\Sigma:\gamma\Rightarrow\gamma$ is contractible (i.e. such that there is a homotopy $h:\Sigma\Rrightarrow \id_\gamma$). This coincides with the identity component of $\hol^2_p$. 
		\end{definition}
		
		\begin{theorem}[Higher Ambrose-Singer]
			The Lie algebra $\mf{hol}^2_p$ of $\hol^2_p$ is given by
			\begin{equation}
				\set[K_{\tra_\gamma(p)}(X,Y,Z)]{\gamma\in \Pi^\thin_1(M),\gamma(0)=x,\,\,X,Y,Z\in T_{\tra(\gamma,p)}P}.
			\end{equation}
		\end{theorem}
		
		\proof Note that $\mf{hol}_p^2$ is spanned by elements of form
		\[
			\left.\frac{\dd}{\dd r}\right|_0\,\tra^2_H(\Sigma_r,p)\tra^2_H(\Sigma_0,p)^{-1},
		\]
		where $\Sigma_r$ is a smooth family of homotopies $\Sigma_r:\gamma\Rightarrow\gamma$. By the higher non-Abelian Stokes theorem \ref{thm:nonabelianstokes2} we have
		\[
		 	= \int_0^1\int_0^1\! \!\dd s\dd t\,\,\Sigma^*K_{0,s,t}(\partial r, \partial s, \partial t) = \int_0^1\int_0^1 \!\! \dd s\dd t\,\,K_{\tra(\Sigma_{0,s,t},p)}(\Sigma_*\partial r,\Sigma_*\partial s,\Sigma_*\partial t).
		\]
		This shows that $\mf{hol}_p^2$ is spanned by the curvature elements. For the converse we note that
		\[
			\left.\frac{\partial}{\partial r}\right|_0\, u_{\rho,\sigma,0}(r,s,\tau +t),
		\]
		as in the proof of theorem \ref{thm:nonabelianstokes2}, is a bigon in $\mf{hol}_p^2$, and is zero if either $s=0$ or $t=0$, hence
		\[
			\left.\frac{\partial^3}{\partial r\partial s \partial t}\right|_{0,0,0}u_{\rho,\sigma,0}(r,s,\tau+t) = K_{\tra(h_{r,s,t},p)}(h_*\partial r,h_*\partial s,h_*\partial t)\in \mf{hol}^2_p.
		\]
		Varying $h:I^3\to M$, and noting that $K$ kills vertical vectors we obtain the inclusion in the other direction. \qed
	
	\addcontentsline{toc}{section}{References}
	\bibliographystyle{alpha}
	\bibliography{2bundles2connection}

\end{document}